\documentclass[AMA,STIX1COL,doublespace]{WileyNJD-v2}
\articletype{Research Article}

\usepackage{amsmath,amsfonts,amsthm}
\usepackage{epsfig}
\usepackage{subfigure}
\usepackage{graphicx}
\usepackage{fancybox}
\usepackage{bm}
\usepackage{verbatim}
\usepackage{cancel}
\usepackage{color,colordvi}
\usepackage{enumerate}

\newcommand{\vm}[1]{\bm{#1}}

\newcommand{\bsym}[1]{\boldsymbol{#1}}

\graphicspath{ {./figs/} }

\begin{document}
	
\title{A Co-rotational Virtual Element Method for 2D Elasticity and Plasticity}
	
\author[1]{Louie L. Yaw$\mbox{}^{*,}$}
	
\authormark{L. L. Yaw}
	
\address[1]{\orgdiv{Engineering Department},
	\orgname{Walla Walla University},
	\orgaddress{College Place, WA 99324, USA}}
\corres{$\mbox{}^*$Louie L. Yaw, Engineering Department,
	Walla Walla University, 100 SW 4th St, College Place, 
	WA 99324, USA \\
	\email{louie.yaw@wallawalla.edu}}

\abstract{
The virtual element method (VEM) allows discretization of the problem domain with polygons in 2D.  The polygons can have an arbitrary number of sides and can be concave or convex.  These features, among others, are attractive for meshing complex geometries.  VEM applied to linear elasticity problems is now well established.  Nonlinear problems involving plasticity and hyperelasticity have also been explored by researchers using VEM.  Clearly, techniques for extending the method to nonlinear problems are attractive.  In this work a novel first order consistent virtual element method is applied within a static co-rotational framework.  To the author’s knowledge this has not appeared before in the literature with virtual elements.  The formulation allows for large displacements and large rotations in a small strain setting.  For some problems avoiding the complexity of finite strains, and alternative stress measures, is warranted.  Furthermore, small strain plasticity is easily incorporated.  The basic method, VEM specific implementation details for co-rotation, and representative benchmark problems are illustrated.  Consequently, this research demonstrates that the co-rotational VEM formulation successfully solves certain classes of nonlinear solid mechanics problems.  The work concludes with a discussion of results for the current formulation and future research directions.}


\maketitle

\section{Introduction}\label{sec:intro}
Since its inception \cite{Beirao,Beirao2} the virtual element method (VEM) has attracted much attention.  Its use in solid mechanics for problems of elasticity \cite{Beirao:2013:VEL,Artioli}, plasticity \cite{TaylorVEMplast,Artioli2,Beirao3}, hyperelasticity \cite{vanHuyssteen,Chi,Wriggers2}, and contact problems \cite{Wriggers3,Aldakheel} has been explored by a variety of researchers.  Clearly, techniques for extending the method to nonlinear problems are attractive.  In this work a novel first order consistent virtual element method is applied within a static co-rotational framework.  To the author's knowledge this has not appeared before in the literature with virtual elements.  The co-rotational formulation, previously applied with finite element \cite{Wempner, Belytschko,Crisfield3, Crisfield4,Jetteur,Rankin1,Rankin2,Rankin3} and meshfree \cite{yawco} simulations, allows for large displacements and large rotations in a small strain setting.  For some problems avoiding the complexity of finite strains, and alternative stress measures, is warranted.  Furthermore, small strain plasticity is easily incorporated.
The structure of the remainder of this work follows.  In section \ref{sec:reviewVEM}, a brief review of virtual elements for 2D elasticity is provided with the reader directed to relevant references for more details.  Furthermore, VEM notation, used herein, is provided.  The VEM review is followed by section \ref{sec:corot}, which provides a derivation and discussion of the co-rotational formulation for linear elastic problems.  This section is essential to identifying how to incorporate VEM into co-rotation and for identifying the key components necessary for nonlinear analysis, namely:  construction of tangent stiffness for each (polygonal) element, calculation of internal force vector for each element, and application of external forces.  By assembling these pieces a typical incremental iterative nonlinear analysis with arc length control is straightforward to implement.  In section \ref{sec:corotplast}, including plasticity in the co-rotational formulation is presented.  Numerical implementation details, particular to this work, are discussed in section \ref{sec:numimplementation}.  Some representative results of numerical simulations are provided in section \ref{sec:numresults}.  Simulation results are compared to benchmark problems or to known theoretical solutions.  Last, main findings and conclusions are provided with thoughts on improvements and future research directions in section \ref{sec:conclusions}.


\section{Review of Virtual Element Method (VEM)}\label{sec:reviewVEM}
The virtual element method (VEM) first appeared in 2013 \cite{Beirao}.  It is a versatile method to solve partial differential equations.  The application of VEM to engineering problems has gained much attention and has much in common with the finite element method (FEM).  A notable difference is that VEM discretizations are not restricted to triangle or quadrilateral polygons; rather, VEM allows for discretizations with arbitrary polygons.  Hence, one can have a collection of triangles, quadrilaterals, pentagons, and so on.  Polygons are permitted to be convex and concave.  This flexibility simplifies mesh construction for problem domains.  The order of polynomial consistency, linear, quadratic, or higher, is allowed if chosen.  Furthermore, VEM is amenable to non-conforming discretizations (see Mengolini \cite{mengolini}).  In this work, attention is restricted to linear order ($k=1$) polynomial interpolation, and this section focuses on solving 2D linear elasticity.    An attempt is made to present minimum implementation details in a  logical order.  The interested reader is directed to the references for details and derivations not provided.  The notation given closely follows references \cite{yawVEM} and \cite{mengolini}, with some exceptions.

Sukumar and Tupek \cite{suku:elastodyn} describe aspects of VEM as follows.  In VEM, ``the basis functions are defined as the solution of a local elliptic partial differential equation", yet they are not ever actually
calculated in construction of the method.  The VEM basis functions are not known and their definition is one of convenience.
This is why they are described as \textit{virtual}, and the finite
element space for VEM as a \emph{virtual element space}.  As is done in the finite element method (FEM), the global conforming space $\boldsymbol{\mathcal{V}}^h$ is constructed by piecing together a local discretization space $\boldsymbol{\mathcal{V}}_k(E)$ (for element E of order $k$). However, differently from FEM, the trial and test functions of $\boldsymbol{\mathcal{V}}_k(E)$
are composed of $k$th order, or less, polynomials, and include nonpolynomial functions as well.
Elliptic polynomial projections of the VEM basis functions are used to construct the stiffness (bilinear form) and forces (linear form) found in a typical weak form.  In linear elasticity problems the degrees of freedom are used to compute the projections with the introduction of no additional approximation error.  With this in hand the stiffness is split into two contributions:  a consistency term associated with the chosen polynomial space and a correction term that provides coercivity (stability or invertibility). The global system of matrices is then put together according to an element by element assembly process analogous to FEM.  One way to imagine VEM is to recognize its similarity to a stabilized hourglass control finite element method~\cite{Flanagan:1981:AUS, Russo} that makes use of convex and or nonconvex polytopes.\cite{Cangiani:2015:HSV}

\subsection{The Continuous 2D Linear Elasticity Problem}
The objective is to solve 2D elasticity problems using VEM.
The standard weak form for elasticity problems is:  Find $\mathbf{u} \in \boldsymbol{\mathcal{V}}$ such that
\begin{equation}\label{E1}
	a(\mathbf{u},\mathbf{v})=L(\mathbf{v}) \quad \forall \mathbf{v}\in \boldsymbol{\mathcal{V}},
\end{equation}
where the bilinear form is
\begin{equation}\label{E2}
	a(\mathbf{u},\mathbf{v})=\int _{\Omega} \! \boldsymbol{\sigma}(\mathbf{u}):\boldsymbol{\varepsilon}(\mathbf{v}) \, \mathrm{d}\Omega,
\end{equation}
and the linear form is
\begin{equation}\label{E3}
	L(\mathbf{v})=\int_{\Omega} \! \mathbf{v}\cdot\mathbf{f} \, \mathrm{d}\Omega+\int_{\partial\Omega_t} \! \mathbf{v}\cdot\bar{\mathbf{t}} ,\ \mathrm{d}\partial\Omega.
\end{equation}\\

\noindent
\textbf{Remarks}
\begin{enumerate}[(i)]
	\item The function space $\boldsymbol{\mathcal{V}}$ is vector-valued with components $v_x$ and $v_y$.  These have zero values on the displacement boundary and are contained in first-order Sobolev space $\mathcal{H}^1(\Omega)$.
	\item The function $\mathbf{u} \in \boldsymbol{\mathcal{V}}$ is a trial solution, $\mathbf{v} \in \boldsymbol{\mathcal{V}}$ is a weight function.
	\item In this work, column vectors, matrices, and Voigt notation are intended in the equations presented, unless it is evident otherwise, such as in \eqref{E2}.
\end{enumerate}

\subsection{Discretization of the problem domain}
\begin{figure}
\centering
\subfigure[]{\epsfig{file =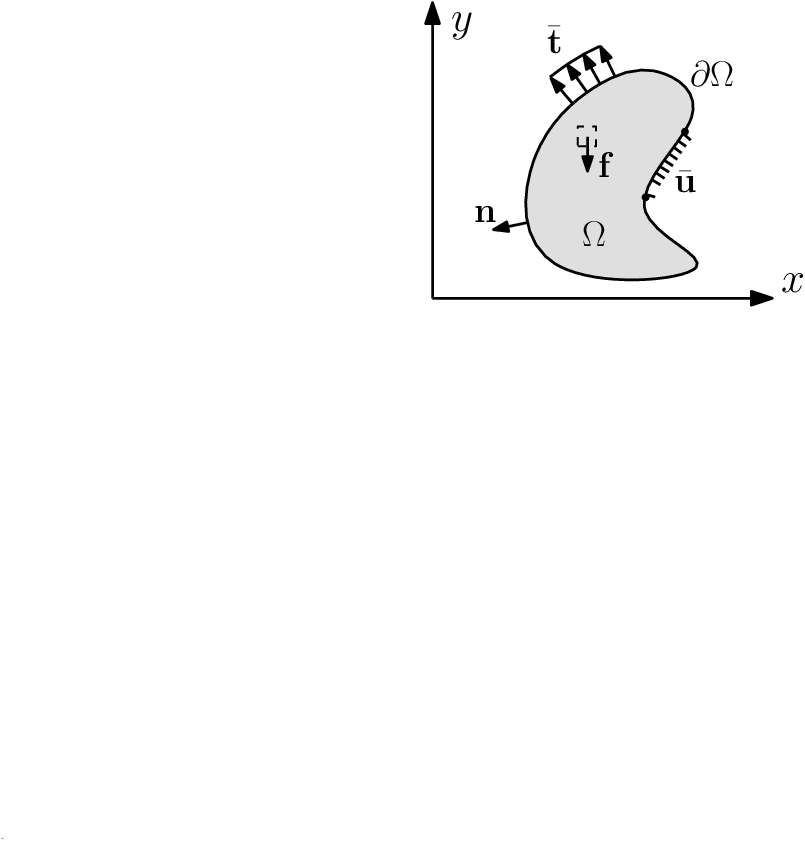,width=0.38\textwidth}}
\hspace*{0.4in}
\subfigure[]{\epsfig{file = 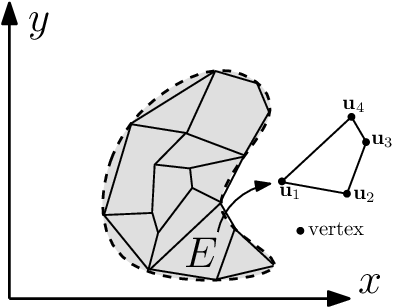,width=0.4\textwidth}}
\caption{2D Solid Domain: (a)~Elasticity problem with boundary conditions, (b)~Virtual element method domain discretization and example polygonal element with vector of nodal displacements labeling each vertex.}
\label{fig1}
\end{figure}

The geometric region shown in Figure~\ref{fig1}a is the geometric domain of interest.  The goal is calculation of field variables such as stress, strain, and displacements.  VEM is employed to accomplish this task by polygonal discretization of the domain, Figure~\ref{fig1}b.  In VEM polygonal elements with an arbitrary number of sides are used and furthermore can be convex or non-convex.  Because of the aforementioned polygonal elements it is necessary to imagine an interpolation space that contains polynomials as well as non-polynomial functions.  This choice allows the sides of polygon elements to compatibly connect along their sides.  As a result, polynomial interpolation functions populate polygon edges, but polynomials plus (possibly) non-polynomial functions reside within the interior of the element.  It so happens that only knowing the polynomial functions along the edges is adequate for success.  Consequently, assuming constant order polynomials for the formulation, a necessary starting point is the choice of polynomial order.  In this work attention is restricted to first order polynomials.
		
A space of scalar-valued polynomials of order equal to $k$ or less on an element $E$ is denoted as $\mathcal{P}_k(E)$.  The space of scalar valued polynomials is extended to a 2D vector space of polynomials with two variables $\boldsymbol{\mathcal{P}}_k\equiv [\mathcal{P}_k]^2$. Hence, a basis for the polynomial space is $\mathbf{P}_k=\{\mathbf{p}_{\alpha}\}_{\alpha=1,...,n_k}$.  Polynomial order $k=1$ is illustrated below.
\begin{equation}\label{pf1}
	\mathbf{P}_1=[\mathbf{p}_1,\ \mathbf{p}_2,\ \mathbf{p}_3,\ \mathbf{p}_4,\ \mathbf{p}_5,\ \mathbf{p}_6]
\end{equation}
or
\begin{equation}\label{pf2}
	\mathbf{P}_1=\left[\left ( \begin{array}{c}
	1\\
	0\\
	\end{array}\right ),
	\ \left ( \begin{array}{c}
	0\\
	1\\
	\end{array}\right ),
	\ \left ( \begin{array}{c}
	-\eta\\
	\xi\\
	\end{array}\right ),
	\ \left ( \begin{array}{c}
	\eta\\
	\xi\\
	\end{array}\right ),
	\ \left ( \begin{array}{c}
	\xi\\
	0\\
	\end{array}\right ),
	\ \left ( \begin{array}{c}
	0\\
	\eta\\
	\end{array}\right )
	\right].
\end{equation}
In equation \eqref{pf2}, the components are constructed by using scaled monomials.  They are defined as
\begin{equation}
	\xi=\left(\frac{x-\bar{x}}{h_E}\right), \quad \eta=\left(\frac{y-\bar{y}}{h_E}\right),
\end{equation}
where $\bar{\mathbf{x}}=(\bar{x},\bar{y})$ is the location of the element $E$ centroid and $h_E$ is the diameter (i.e., diameter of smallest circle that encloses all element vertices).\\

\noindent
\textbf{Remarks}
\begin{enumerate}[(i)]
	\item A choice of higher order polynomials is possible (see \cite{mengolini}).
	\item Typically the analyst chooses the location of degrees of freedom.  For example, the polygon vertices are chosen in this work so that each vertex includes an $x$ and $y$ displacement degree of freedom.  Consequently,  an element edge has two vertices (points), which establish a line, and a first order polynomial is consistent with this circumstance.  The $x$ values, for example, can be linearly interpolated along the element edges.
	\item For 2D elasticity $[\mathcal{P}_k]^2$ denotes a 2D vector polynomial with components in two variables.
	\item A polynomial base has a number of terms equal to $n_k=(k+1)(k+2)$.  For $k=1$, $n_k=6$.  The polynomial base $\mathbf{P}_1$ contains 6 terms.
	\item The first three monomials $\mathbf{p}_1, \mathbf{p}_2, \mathbf{p}_3$ of \eqref{pf2} provide terms associated with rigid body motion.
	\item For $\alpha=1,2,3$ of the polynomial base the infinitesimal strain equals zero.  Later, this has consequences for the construction of the $\bar{\mathbf{B}}$ matrix of \eqref{eqbbar}, where the first three rows are populated with \eqref{p16}, otherwise the first three rows would all be zero (see \cite{yawVEM, mengolini} for more details).  The $\bar{\mathbf{B}}$ matrix mentioned here and defined in subsection \ref{stiffsect} is not to be confused with the $\bar{\mathbf{B}}$ matrix used in the finite element method for problems of incompressibility \cite{Hughes}.
	\item Henceforth, $n_v$ denotes the number of vertices for a given polygon.  For a first order consistent formulation, the number of degrees of freedom in one spatial direction, $n_d$, is equivalent to $n_v$.
\end{enumerate}	

\subsection{Element Stiffness}\label{stiffsect}
As with finite element analysis, in VEM an element stiffness matrix is constructed for each polygonal element.  The element stiffness in VEM is composed of two parts: a consistency part and a stability part.
\begin{equation}
	\mathbf{k}_E=\mathbf{k}_E^c+\mathbf{k}_E^s \quad \quad [2n_v \times 2n_v]
\end{equation}
The consistency part is expressed as
\begin{equation}
	\mathbf{k}_E^c=tA_E\mathbf{B}^T\mathbf{C}\mathbf{B},
\end{equation}
where $\mathbf{C}$ is the elastic modular matrix, $\mathbf{B}$ is the strain displacement matrix, $A_E$ is the area of polygon element $E$, and $t$ is the thickness of the polygonal element.  For VEM the strain displacement matrix is written
\begin{equation}
	\mathbf{B}=\bsym{\varepsilon}(\mathbf{P}_1)\tilde{\mathbf{\Pi}}. \quad \quad [3 \times 2n_v]
\end{equation}
In the above, the strain operator, $\bsym{\varepsilon}$, acts on the polynomial basis in the following way:
\begin{equation}
	\bsym{\varepsilon}(\mathbf{P}_1)=\bsym{\varepsilon}([\mathbf{p}_1 \ \mathbf{p}_1 \ \hdots \ \mathbf{p}_{n_k}])
\end{equation}
Specifically, the strain operator acts on the polynomial base functions and creates the matrix,
\begin{equation}\label{ep8}
	\begin{split}
		\bsym{\varepsilon}\left[ \begin{array}{cccc}
		\mathbf{p}_1 & \mathbf{p}_2 & ... & \mathbf{p}_{n_k} 
		\end{array}\right]=\left[ \begin{array}{cccc}
		\partial_{x}p_{1,1} & \partial_{x}p_{2,1} & ... & \partial_{x}p_{n_k,1} \\
		\partial_{y}p_{1,2} & \partial_{y}p_{2,2} & ... & \partial_{y}p_{n_k,2} \\
		\partial_{y}p_{1,1}+\partial_{x}p_{1,2} & \partial_{y}p_{2,1}+\partial_{x}p_{2,2} & ... & \partial_{y}p_{n_k,1}+\partial_{x}p_{n_k,2} 
		\end{array}\right],
	\end{split}
\end{equation}
where in \eqref{ep8}, $p_{i,j}$ is used to represent the $j$th (1st or 2nd) component of polynomial vector $\mathbf{p}_i$.  In 2D the $\mathbf{p}_i$ are expressed in equations \eqref{pf1} and \eqref{pf2}.  The engineering strains, themselves, are obtained using the element vertex displacements, $\mathbf{u}^E$, as shown below:
\begin{equation}\label{ep7}
	\begin{split}
		\bsym{\varepsilon}(\mathbf{v}^h)\approx\bsym{\varepsilon}(\Pi(\mathbf{v}^h))&=\bsym{\varepsilon}\left(\left[ \begin{array}{cccc}
		\mathbf{p}_1 & \mathbf{p}_2 & ... & \mathbf{p}_{n_k} 
		\end{array}\right]\tilde{\bsym{\Pi}}\mathbf{u}^E\right)\\
		&=\bsym{\varepsilon}\left[ \begin{array}{cccc}
		\mathbf{p}_1 & \mathbf{p}_2 & ... & \mathbf{p}_{n_k} 
		\end{array}\right]\tilde{\bsym{\Pi}}\mathbf{u}^E\\
		&=\mathbf{B}\mathbf{u}^E. \quad \quad [3 \times 1]
	\end{split}
\end{equation} 
\noindent
\textbf{Remarks}
\begin{enumerate}[(i)]
	\item Voigt notation is used to organize the strains.  Two-dimensional engineering strains $\varepsilon_x$, $\varepsilon_y$, $\gamma_{xy}=2\varepsilon_{xy}$ populate the strain vector $\bsym{\varepsilon}(\mathbf{v}^h)$.  The strain displacement operator $\mathbf{B}$ appears in~\eqref{ep7} and is $3 \times 2n_v$ in size.
	\item The displacement vector is organized in the standard way, $\mathbf{u}^E=[u^1_x \ \ u^1_y \ \ u^2_x \ \ u^2_y \cdots u^{n_v}_x \ \ u^{n_v}_y]^T$.
	\item Strains are constant within an individual polygon element.
	\item The expression $\Pi(\mathbf{v}^h)$ uses the projection operator $\Pi$ and implies the projection of the virtual element approximation, $\mathbf{v}^h$, onto the polynomial space.
	\item The strain operator $\bsym{\varepsilon}$ is context specific.  If it operates on a single vector, like $\mathbf{p}_1$ or $\mathbf{v}^h$, in 2D with two components, then the result is $3 \times 1$.  If it operates on a group of vectors, like $\mathbf{P}_1$, which is happening in \eqref{ep8}, then the result is $3 \times n_k$ where, for the 1st order polynomial base, $n_k=6$.
\end{enumerate}
The projector matrix, $\tilde{\mathbf{\Pi}}$, is constructed as
\begin{equation}\label{eqbbar}
	\tilde{\mathbf{\Pi}}=\mathbf{G}^{-1}\bar{\mathbf{B}},  \quad \quad [n_k \times 2n_v]
\end{equation}
where
\begin{equation}
	\mathbf{G}=\bar{\mathbf{B}}\mathbf{D}.  \quad \quad [n_k \times n_k]
\end{equation}
The $\bar{\mathbf{B}}$ matrix is assembled by forming the $\tilde{\mathbf{B}}$ matrix and then replacing the first three rows with the $\breve{\mathbf{B}}$ matrix.  The $\tilde{\mathbf{B}}$ matrix is formed, for rows $\alpha=1$ to $n_k$ two column at a time, for $j$ values ranging over the number of vertices $1$ to $n_v$, by using the following expressions:
\begin{equation}\label{ne11}
	\begin{split}
		\tilde{B}_{\alpha (2j-1)}&=\left[ \begin{array}{c}
		1\\
		0
		\end{array}\right]\cdot \left[ \begin{array}{cc}
		\sigma_x(\mathbf{p}_{\alpha}) & \sigma_{xy}(\mathbf{p}_{\alpha})\\
		\sigma_{xy}(\mathbf{p}_{\alpha}) &\sigma_y(\mathbf{p}_{\alpha})
		\end{array}\right]
		\left(\frac{|e_{j-1}|}{2}
		\left[ \begin{array}{c}
		n_{e1}\\
		n_{e2}
		\end{array}\right]_{j-1}+\frac{|e_{j}|}{2}
		\left[ \begin{array}{c}
		n_{e1}\\
		n_{e2}
		\end{array}\right]_j\right)\\ 
		\text{\normalfont and}&\\
		\tilde{B}_{\alpha (2j)}&=\left[ \begin{array}{c}
		0\\
		1
		\end{array}\right]\cdot \left[ \begin{array}{cc}
		\sigma_x(\mathbf{p}_{\alpha}) & \sigma_{xy}(\mathbf{p}_{\alpha})\\
		\sigma_{xy}(\mathbf{p}_{\alpha}) &\sigma_y(\mathbf{p}_{\alpha})
		\end{array}\right]
		\left(\frac{|e_{j-1}|}{2}
		\left[ \begin{array}{c}
		n_{e1}\\
		n_{e2}
		\end{array}\right]_{j-1}+\frac{|e_{j}|}{2}
		\left[ \begin{array}{c}
		n_{e1}\\
		n_{e2}
		\end{array}\right]_j\right).
	\end{split}
\end{equation}
	\begin{figure}
  \centering
  \epsfig{file=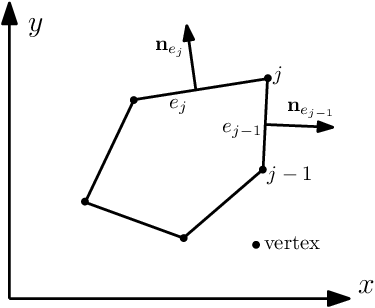,width=.4\textwidth}
  \caption{Single five sided element:  edges, normals, and nodes labeled.}\label{fige}
  \end{figure}
For equation \eqref{ne11}, with reference to Figure \ref{fige}, edge length $j$ is written as $|e_j|$ and the $j$th edge outward unit normal vector components are $[n_{e1} \ n_{e2}]^T_j$.  For vertex $j = 1$, edge length $|e_{j-1}|$ is the length from vertex $j = n_v$ to $j = 1$, where the polygon element has $n_v$ vertices, and normal $n_{e_{j-1}}$ is the outward unit normal to the edge from $j = n_v$ to $j = 1$.  Also, note that the right side of \eqref{ne11} yields a vector that is subsequently dotted with $[1 \ \ 0]^T$ or $[0 \ \ 1]^T$, as indicated.	
		
The engineering stress terms, $\bsym{\sigma}(\mathbf{p}_{\alpha})$, are found by matrix multiplication
\begin{equation}\label{ne8}
	\begin{split}
		\bsym{\sigma}(\mathbf{p}_{\alpha})&=\mathbf{C}\boldsymbol{\varepsilon}(\mathbf{p}_{\alpha})\\
		\left[ \begin{array}{c}
		\sigma_x(\mathbf{p}_{\alpha})\\
		\sigma_y(\mathbf{p}_{\alpha})\\
		\sigma_{xy}(\mathbf{p}_{\alpha})
		\end{array}\right]&=\mathbf{C}
		\left[ \begin{array}{c}
		\varepsilon_x(\mathbf{p}_{\alpha})\\
		\varepsilon_y(\mathbf{p}_{\alpha})\\
		\gamma_{xy}(\mathbf{p}_{\alpha})
		\end{array}\right].
	\end{split}
\end{equation}
In \eqref{ne8}, the appropriate $\mathbf{C}$ matrix for plane stress or plain strain is used (see equations \eqref{s2} and \eqref{s3}).  The result of \eqref{ne8} is used to form the stress matrix
\begin{equation}\label{ne9}
	\left[\boldsymbol{\sigma}(\mathbf{p}_{\alpha})\right]=\left[ \begin{array}{cc}
		\sigma_x(\mathbf{p}_{\alpha}) & \sigma_{xy}(\mathbf{p}_{\alpha})\\
		\sigma_{xy}(\mathbf{p}_{\alpha}) &\sigma_y(\mathbf{p}_{\alpha})
		\end{array}\right].
\end{equation}
Then \eqref{ne9} is used in \eqref{ne11}.  For plane stress
\begin{equation}\label{s2}
	\mathbf{C}=\frac{E_Y}{1-\nu^2}\begin{bmatrix}
		\phantom{0}1 & \phantom{00}\nu & \phantom{00}0\\
		\phantom{0}\nu & \phantom{00}1 & \phantom{00}0\\
		\phantom{0}0 & \phantom{00}0 & \phantom{00}\frac{1-\nu}{2}
		\end{bmatrix},
\end{equation}
and for plain strain
\begin{equation}\label{s3}
	\mathbf{C}=\frac{E_Y}{(1+\nu)(1-2\nu)}\begin{bmatrix}
		1-\nu & \nu & 0\\
		\nu & 1-\nu & 0\\
		0 & 0 & \frac{1-2\nu}{2}
		\end{bmatrix}\color[rgb]{0,0,0},
\end{equation}
where $E_Y$ is Young's modulus of elasticity, and $\nu$ is Poisson's ratio.

The calculation of \eqref{ne8} describes `stress' terms necessary for the construction of entries in the $\tilde{\mathbf{B}}$ matrix.  However, to calculate actual engineering stresses during post-processing for elasticity problems, the stress calculation makes use of \eqref{ep7} and the appropriate modular matrix $\mathbf{C}$.  The result is
\begin{equation}\label{stress1}
	\bsym{\sigma}(\mathbf{v}^h)\approx \mathbf{C} \bsym{\varepsilon}(\Pi(\mathbf{v}^h))
	=\mathbf{C}\mathbf{B}\mathbf{u}^E.
\end{equation}

The $\breve{\mathbf{B}}$ matrix is calculated as 
\begin{equation}\label{p16}
	\breve{\mathbf{B}}= \left[ \begin{array}{cccc}
		\breve{\mathbf{b}}_1 & \breve{\mathbf{b}}_2& ... & \breve{\mathbf{b}}_{2n_v} \\
		\end{array}\right],   \quad \quad [3 \times 2n_v]
\end{equation}
where
\begin{equation}\label{p15}
	\breve{\mathbf{b}}_I=\left [ \begin{array}{c}
		\frac{1}{n_v}\sum\limits^{2n_v}_{i=1}\delta_{iI}dof_i(\mathbf{p}_1)\\[0.3cm]
		\frac{1}{n_v}\sum\limits^{2n_v}_{i=1}\delta_{iI}dof_i(\mathbf{p}_2)\\[0.3cm]
		\frac{1}{n_v}\sum\limits^{2n_v}_{i=1}\delta_{iI}dof_i(\mathbf{p}_3)
		\end{array}\right ], \quad \quad \delta_{iI}=\left\{ \begin{array}{c} 
		0 \ \text{for} \ i\neq I\\
		1 \ \text{for} \ i=I.
		\end{array}
		\right.
\end{equation}
Note that, $dof_i(\mathbf{p}_j)$ is found as follows:  first, calculate polynomial vector $\mathbf{p}_j$ evaluated at the vertex coordinates associated with dof $i$, second, take the component of the vector that is directed along dof $i$.  Essentially, $dof_i$ is an operator that extracts the value of its argument, at dof $i$, in the direction of dof $i$, and the dofs (degrees of freedom) range from $1$ to $2n_v$.

The $\mathbf{D}$ matrix is found by computing polynomial vector values at the coordinates associated with the degrees of freedom of polygon $E$.  The result is
\begin{equation}
	\mathbf{D}=\left[ \begin{array}{cccc}
		dof_1(\mathbf{p}_{1}) & dof_1(\mathbf{p}_{2})& \cdots & dof_1(\mathbf{p}_{n_k})\\
		dof_2(\mathbf{p}_{1}) & dof_2(\mathbf{p}_{2})& \cdots & dof_2(\mathbf{p}_{n_k})\\
		\vdots & \vdots & \ddots & \vdots\\
		dof_{2n_v}\mathbf{p}_{1}) & dof_{2n_v}(\mathbf{p}_{2})& \cdots & dof_{2n_v}(\mathbf{p}_{n_k})
		\end{array}\right].
\end{equation}

Researchers prescribe the stability part of the element stiffness in a variety of forms.  One form suggested by Menglolini~et~al. \cite{mengolini} is
\begin{equation}
	\mathbf{k}_E^s=\tau^h \ \textrm{tr}(\mathbf{k}_E^c)(\mathbf{I}-\bsym{\Pi})^T(\mathbf{I}-\bsym{\Pi}),
\end{equation}
where $\text{\normalfont tr}$~denotes the trace operator, and $\tau^h=1/2$ is suggested for linear elasticity. The identity matrix $\mathbf{I}$ is $2n_v \times 2n_v$.  The projection matrix, $\bsym{\Pi}$, is determined as
\begin{equation}\label{stab7}
	\bsym{\Pi}=\mathbf{D}\tilde{\bsym{\Pi}}.  \quad \quad [2n_v \times 2n_v]
\end{equation}
Another form of the stability part of the element stiffness, suggested by Sukumar and Tupek \cite{suku:elastodyn}, is expressed as
\begin{equation}
	\mathbf{k}_E^s=(\mathbf{I}-\bsym{\Pi})^T\mathbf{S}^d_E(\mathbf{I}-\bsym{\Pi}).
\end{equation}
where a $2n_v \times 2n_v$ diagonal matrix, $\mathbf{S}^d_E$, is scaled as needed.  The diagonal terms of the matrix are prescribed as: $(\mathbf{S}^d_E)_{ii}=\mathsf{max}(\alpha_0 \ \text{\normalfont tr}(\mathbf{C})/m, (\mathbf{k}^c_E)_{ii})$, where $m=3$ in 2D, the modular matrix $\mathbf{C}$ in 2D is as appropriate for plane strain or plane stress, and $\alpha_0=1$ because the scaled monomials of VEM create polygon elements with diameters on the order of 1.

Both forms of the stability part of the element stiffness matrix above are found to be effective, with negligible difference between the two, for compressible materials in linear elasticity problems.  They are also effective for plane stress plasticity problems.  Though not included in the numerical examples of this work, for researchers interested in problems of near incompressibility, the author finds an approach like \cite{park} to be very effective. 

\subsection{Application of External Forces}
Tractions and body forces cause external forces as given in equation \eqref{E3}.  It is possible to apply point loads at individual nodes (vertices) of polygon elements also.  Application of forces is accomplished in the same manner as in standard finite element methods.  The three types of forces are included in the linear form for a single element.
\begin{equation}\label{ext1}
	L_E(\mathbf{v}^h)=\int_E \! \mathbf{v}^h\cdot\mathbf{f} \, \mathrm{d}E+\int_{\partial E\cap\Omega_t} \! \mathbf{v}^h\cdot\bar{\mathbf{t}} \, \mathrm{d}\partial E+\sum\limits_{i=1}\mathbf{v}^h(\mathbf{x}_i)\cdot\mathbf{F}_i. 
\end{equation}
For more discussion the reader is directed to Mengolini et al. \cite{mengolini}.  In the examples, shown later in the results section, only point loads are used.  Often the point load is evenly distributed over a small region or line of vertices along an edge.  This is sufficient for the example problems provided.  From point loads an external force vector is assembled, which is used to solve for the nodal displacements.  In a nonlinear analysis, the global external forces are used in a Newton-Raphson scheme to enforce equilibrium with the global internal force vector \eqref{eif2}.  The word `global' is used to create a distinction from quantities in local coordinates.  This distinction becomes important in a co-rotational formulation.

\section{Co-rotational VEM - Elasticity}\label{sec:corot}
\begin{figure}
\centering
\epsfig{file=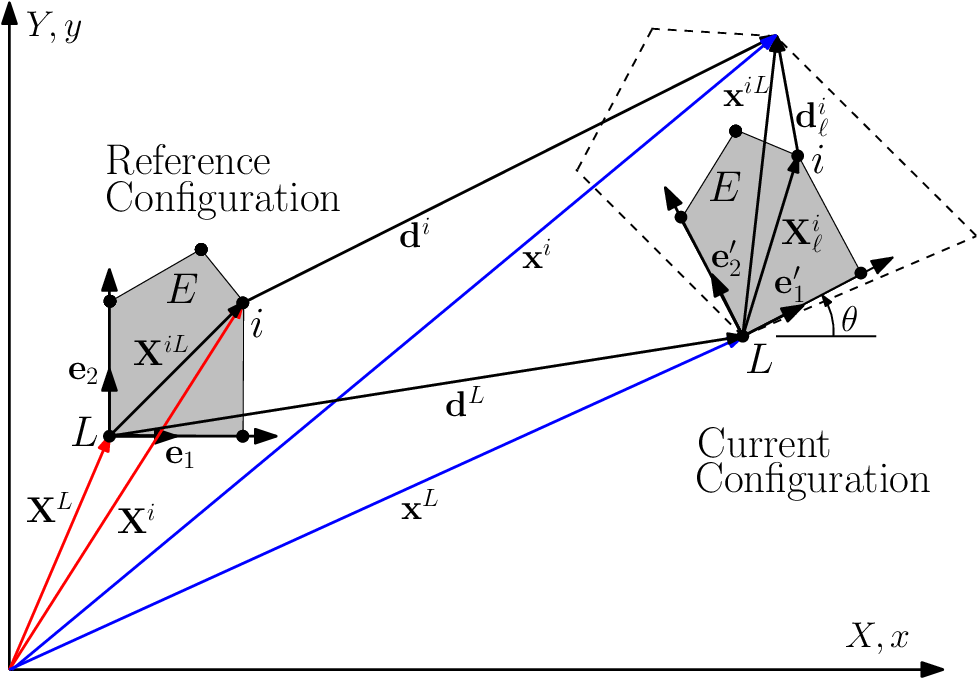,width=.65\textwidth}
\caption{Reference and current configurations in co-rotational formulation.}\label{RCconfig}
\end{figure}

To extend VEM to include large displacements and large rotations a co-rotational formulation is employed.  In this section, the discussion and derivations closely follow \cite{Crisfield3} and \cite{yawco}, modified as needed for VEM related considerations.  As a polygonal element rotates and translates during deformations a local co-rotating frame is attached to node $L$ (typically node 1) of the element (see Figure~\ref{RCconfig}).  The angle of rotation, $\theta$, of the co-rotating frame is determined in the current configuration.  The local co-rotating frame and the global axes start parallel to each other in the reference configuration.  Translations are found by knowing reference coordinates and current coordinates for the element nodes.  With these items in hand rigid body translations and rotations are removed from the current displacements.  In the local co-rotating coordinate system only strain causing deformations remain.  As a result, in the local co-rotating frame, small strain elasticity (or plasticity) elements are easily incorporated.  The key ingredients of a co-rotational formulation are: (i) the relationship between local and global variables, (ii) the angle of the local co-rotating frame, and (iii) a variationally consistent tangent stiffness matrix.

\subsection{Relationship between global and local variables}
Figure~\ref{RCconfig} is used to illustrate the relations between global displacements and local strains.  For a particular element, node $L$ and the remaining element nodes are shown.  The displacement of the overall structure causes the polygon element to translate, rotate, and deform to the current configuration.  For node $i$, the initial and current global nodal position vectors are denoted as $\mathbf{X}^i$ and $\mathbf{x}^i$.  A subscript $\ell$ is attached to vectors in local coordinates, such as local nodal vectors, $\mathbf{X}^i_{\ell}$, and local nodal displacement vectors, $\mathbf{d}^i_{\ell}$ (some exceptions are stress, $\bsym{\sigma}$, and strain-displacement matrices, $\mathbf{B}$, understood as being in the local coordinates of the element's co-rotational frame).  Current global nodal displacements vectors are written as $\mathbf{d}^i$.  The current difference between nodes $i$ and $L$ is written as $\mathbf{x}^{iL}=\mathbf{x}^i-\mathbf{x}^L=\mathbf{X}^{iL}+\mathbf{d}^i-\mathbf{d}^L$.  An orthogonal orientation matrix, $\mathbf{Q}=[\mathbf{e}_1' \ \mathbf{e}_2']$, for the local co-rotating frame is written in terms of the local unit basis vectors $\mathbf{e}_1'$ and $\mathbf{e}_2'$.  Global vector quantities are transformed to local vector quantities by using $\mathbf{Q}^T$.  The current angle, $\theta$, of the co-rotating frame is used to express the local basis vectors with global components as

\begin{equation}\label{Rgd2}
\mathbf{e}'_1=\left[\begin{array}{c}
\cos{\theta} \\
\sin{\theta} 
\end{array} \right] \quad \mathbf{e}'_2=\left[\begin{array}{c}
-\sin{\theta} \\
\cos{\theta} 
\end{array}\right].
\end{equation}

Using the above, the local nodal displacement vectors for each node $i$ are
\begin{equation}\label{Rgd1}
\mathbf{d}^i_{\ell}=\mathbf{Q}^T\mathbf{x}^{iL}-\mathbf{X}^i_{\ell},
\end{equation}
where $\mathbf{X}^i_{\ell}=\mathbf{Q}^T\mathbf{X}^{iL}=\mathbf{Q}^T\left(\mathbf{X}^i-\mathbf{X}^L\right)$.

\subsection{Angle of the co-rotating frame}\label{sec:angle}
The angle of rotation, due to current local nodal displacements, is found by requiring zero spin at the element centroid (see Jetteur and Cescotto \cite{Jetteur}).  (Alternatively, this could be accomplished by a polar decomposition approach.)

\begin{equation}\label{angle1}
\Omega_{\ell}=\frac{\partial u_{1\ell}}{\partial Y_{\ell}}-\frac{\partial u_{2\ell}}{\partial X_{\ell}}=0.
\end{equation}
The local displacement field in the direction of degree of freedom $j$ is written  as
\begin{equation}\label{angle2}
 u_{j\ell}=\boldsymbol{\phi}^T\mathbf{d}_{j\ell}, \quad (j=1,2)
\end{equation}
where $\mathbf{d}_{j\ell}$ is the vector of local element displacements for degree of freedom $j$ and $\bsym{\phi}$ is the vector of element basis functions.  The basis functions for VEM are not known, yet the approximation of their derivatives, needed to calculate $\mathbf{a}_\ell$ (see \eqref{angle3}b), are present within the strain displacement matrix, $\mathbf{B}$.  These derivatives, so obtained, prove to give satisfactory results for the co-rotational formulation.
Substituting \eqref{angle2} into \eqref{angle1} yields
\begin{subequations}\label{angle3}
\begin{align}
\Omega_{\ell}&=\left(\frac{\partial \boldsymbol{\phi}}{\partial Y_{\ell}}\right)^T \mathbf{d}_{1\ell} - \left(\frac{\partial \boldsymbol{\phi}}{\partial X_{\ell}}\right)^T \mathbf{d}_{2\ell}=\mathbf{a}^T_{\ell}\mathbf{d}_{\ell}, \\
\intertext{where}
\mathbf{a}_{\ell}&=\left[\begin{array}{c}
\frac{\partial \boldsymbol{\phi}_1}{\partial Y_{\ell}} \\ \noalign{\medskip}
-\frac{\partial \boldsymbol{\phi}_1}{\partial X_{\ell}}\\ \noalign{\medskip}
\vdots \\ \noalign{\medskip}
\frac{\partial \boldsymbol{\phi}_n}{\partial Y_{\ell}}\\ \noalign{\medskip}
-\frac{\partial \boldsymbol{\phi}_n}{\partial X_{\ell}}
\end{array} \right] \quad \mathrm{and} \quad \mathbf{d}_{\ell}=\left[\begin{array}{c}
d^{1}_{1\ell} \\ 
 d^{1}_{2\ell}\\ 
d^{2}_{1\ell} \\ 
 d^{2}_{2\ell}\\ 
 \vdots \\ 
d^{n}_{1\ell} \\ 
d^{n}_{2\ell}
\end{array} \right].
\end{align}
\end{subequations}
In equation \eqref{angle3}b, local node $s$ displacement in the direction of dof $r$ ($1=x$ and $2=y$ direction) is written as $d^s_{r\ell}$.
For the above expressions $n=n_v$ is the number of vertices in the polygon element.  Substitute \eqref{Rgd1} into \eqref{angle3}a to obtain
\begin{equation}\label{angle6}
\Omega_{\ell}=\sum(\mathbf{a}^i_{\ell})^T(\mathbf{Q}^T\mathbf{x}^{iL})
-\sum(\mathbf{a}^i_{\ell})^T(\mathbf{X}^i_{\ell})=0.
\end{equation}
The last term in \eqref{angle6} is zero.  The expanded first term is
\begin{equation*}\label{angle7}
\Omega_{\ell}=\sum(\mathbf{a}^i_{\ell})^T\left(\cos{\theta}\left( \begin{array}{c}
x^{iL} \\ 
y^{iL}
\end{array} \right)+\sin{\theta}\left( \begin{array}{c}
y^{iL} \\ 
-x^{iL}
\end{array} \right) \right)
=0,
\end{equation*}
or
\begin{equation}\label{angle8}
\Omega_{\ell}=a\sin{\theta}+b\cos{\theta}=0,
\quad a=\sum(\mathbf{a}^i_{\ell})^T\left( \begin{array}{c}
y^{iL} \\ 
-x^{iL}
\end{array} \right), \quad b=\sum(\mathbf{a}^i_{\ell})^T\left( \begin{array}{c}
x^{iL} \\ 
y^{iL}
\end{array} \right).
\end{equation}
The results of \eqref{angle8} are also expressed as:
\begin{equation}\label{angle10}
a=\mathbf{c}^T\hat{\mathbf{x}},
\end{equation}
where
\begin{equation}\label{angle11}
\mathbf{c}=\left[\begin{array}{ccccccc}
0 & -1 & 0 & 0 & \hdots & 0 & 0 \\ 
1 & 0 & 0 & 0 &   & 0 & 0 \\ 
0 & 0 & 0 & -1 &  &  &  \vdots\\ 
0 & 0 & 1 & 0 &  &  &  \\ 
\vdots &  &  &  & \ddots &  &  \\ 
0 & 0 &  &  &  & 0 & -1 \\ 
0 & 0 & \hdots &  &  & 1 & 0
\end{array}  \right]\mathbf{a}_{\ell}, \quad \hat{\mathbf{x}}=\left[\begin{array}{c}
 x^{1L} \\y^{1L} \\x^{2L}  \\y^{2L}  \\ \vdots  \\ x^{nL} \\y^{nL}  
\end{array}\right]
\end{equation}
and $b=\mathbf{a}^T_{\ell}\hat{\mathbf{x}}$.  The matrix $\mathbf{c}$ in \eqref{angle11} is $2n \times 2n$ depending on the number $n=n_v$ of polygonal vertices.  The column vector $\hat{\mathbf{x}}$ is $2n \times 1$.  With the above and $\eqref{angle8}$ in hand, the angle of rotation is
\begin{equation}\label{angle12}
\theta=\tan^{-1}\left( \frac{-b}{a} \right).
\end{equation}

\subsection{Derivation of the tangent stiffness matrix}
The tangent stiffness matrix is derived by considering the local internal force vector, $\mathbf{q}_{\ell}$, for element $E$.  It is written as
\begin{equation}\label{kt1}
\mathbf{q}_{\ell}=\int_{\Omega} \mathbf{B}^T\boldsymbol{\sigma} \ dV=\mathbf{k}_{t\ell}\mathbf{d}_{\ell},
\end{equation}
where the local strain-displacement matrix is $\mathbf{B}$, the local engineering stresses are $\boldsymbol{\sigma}$, and the local material tangent stiffness matrix is represented as $\mathbf{k}_{t\ell}$.

Next, by way of some function, $f$, local nodal displacements, $\mathbf{d}_{\ell}$, relate to global nodal displacements, $\mathbf{d}$, and the rotation basis vectors, i.e.,
\begin{equation}\label{kt2}
\mathbf{d}_{\ell}=f(\mathbf{d},\mathbf{e}'_1,\mathbf{e}'_2).
\end{equation}
The variation of \eqref{kt2} produces
\begin{equation}\label{kt3}
\delta\mathbf{d}_{\ell}=\mathbf{T}\delta\mathbf{d},
\end{equation}
where $\mathbf{T}$ is a transformation matrix, which is to be determined.
Recognizing the equivalence of virtual work at the local and global level yields
\begin{equation}\label{kt4oa}
	(\delta\mathbf{d}_{\ell})^T\mathbf{q}_{\ell}=(\delta\mathbf{d})^T\mathbf{q}.
\end{equation}
By making use of \eqref{kt1}, \eqref{kt3} and \eqref{kt4oa}, global internal forces in terms of local internal forces are found, giving
\begin{equation}\label{kt4}
\mathbf{q}=\mathbf{T}^T\mathbf{q}_{\ell}=\mathbf{T}^T\mathbf{k}_{t\ell}\mathbf{d}_{\ell}.
\end{equation}
Then taking the variation of \eqref{kt4}, the global stiffness matrix is found next,
\begin{equation}\label{kt5}
\delta\mathbf{q}=\mathbf{T}^T\delta\mathbf{q}_{\ell}+\delta\mathbf{T}^T\mathbf{q}_{\ell}
=\mathbf{T}^T\mathbf{k}_{t\ell}\delta\mathbf{d}_{\ell}+\mathbf{k}_{t\sigma}\delta\mathbf{d}
=\mathbf{T}^T\mathbf{k}_{t\ell}\mathbf{T}\delta\mathbf{d}+\mathbf{k}_{t\sigma}\delta\mathbf{d},
\end{equation}
where $\delta\mathbf{T}^T\mathbf{q}_{\ell}$ is denoted by $\mathbf{k}_{t\sigma}\delta\mathbf{d}$.  The initial stiffness matrix is $\mathbf{k}_{t\sigma}$, the local material tangent stiffness is $\mathbf{k}_{t\ell}$ (for the case of linear elasticity or possibly constructed by considering inelastic material behavior).  By making use of \eqref{kt3} the last equality in \eqref{kt5} is reached.  Equation \eqref{kt5} yields
\begin{equation}\label{kt6}
\delta\mathbf{q}=\left[\mathbf{T}^T\mathbf{k}_{t\ell}\mathbf{T}+\mathbf{k}_{t\sigma}\right]\delta\mathbf{d}
=\mathbf{k}_T\delta\mathbf{d},
\end{equation}
where the global tangent stiffness matrix for an element is $\mathbf{k}_T$.  It is important to deliberately point out that the VEM small strain stiffness $\mathbf{k}_E$ \emph{is} the material stiffness that is inserted in place of $\mathbf{k}_{t\ell}$ in equation \eqref{kt6}.

To determine the transformation matrix in \eqref{kt6}, it is necessary to take the variation of \eqref{Rgd1}, which yields
\begin{equation}\label{kt7}
\delta\mathbf{d}^i_{\ell}=\mathbf{Q}^T\delta\mathbf{x}^{iL}+\delta\mathbf{Q}^T\mathbf{x}^{iL}.
\end{equation}
Recall from Figure \ref{RCconfig}, that
\begin{equation}\label{kt8}
\mathbf{x}^{iL}=\mathbf{X}^{iL}+\mathbf{d}^i-\mathbf{d}^L=\mathbf{X}^{iL}+\mathbf{d}^{iL}.
\end{equation}
From the variation of \eqref{kt8}, obtain
\begin{equation}\label{kt9}
\delta\mathbf{x}^{iL}=\delta\mathbf{X}^{iL}+\delta\mathbf{d}^{iL}=\delta\mathbf{d}^{iL},
\end{equation}
where the last equality occurs because $\delta\mathbf{X}^{iL}$ is zero.  Inserting \eqref{kt9} into \eqref{kt7} gives
\begin{equation}\label{kt10}
\delta\mathbf{d}^i_{\ell}=\mathbf{Q}^T\delta\mathbf{d}^{iL}+\delta\mathbf{Q}^T\mathbf{x}^{iL}.
\end{equation}
The variation of $\mathbf{Q}^T$ yields
\begin{equation}\label{kt11}
\delta\mathbf{Q}^T=\delta\left[\begin{array}{cc}
\mathbf{e}'_1 & \mathbf{e}'_2
\end{array} \right]^T=\left[\begin{array}{cc}
-\sin{\theta} &-\cos{\theta} \\ 
\cos{\theta} &-\sin{\theta}
\end{array} \right]^T\delta\theta.
\end{equation}
Consequently, (using $s=\sin{\theta}$ and $c=\cos{\theta}$)
\begin{equation}\label{kt12}
\delta\mathbf{Q}^T\mathbf{x}^{iL}=\left[\begin{array}{cc}
 -s&c  \\ 
 -c&-s 
\end{array} \right]
\left(\begin{array}{c}
x^{iL} \\ 
y^{iL}
\end{array} \right)\delta\theta=\left[\begin{array}{c}
-sx^{iL}+cy^{iL} \\ 
-cx^{iL}-sy^{iL}
\end{array} \right]\delta\theta=\mathbf{Q}^T\left(\begin{array}{c}
y^{iL} \\ 
-x^{iL}
\end{array} \right)\delta\theta.
\end{equation}
Next, inserting \eqref{kt12} into \eqref{kt10} gives
\begin{equation}\label{kt13}
\delta\mathbf{d}^i_{\ell}=\mathbf{Q}^T\delta\mathbf{d}^{iL}+\mathbf{Q}^T\left(\begin{array}{c}
y^{iL} \\ 
-x^{iL}
\end{array} \right)\delta\theta.
\end{equation}
By adding $\mathbf{Q}^T\delta\mathbf{d}^L$ to \eqref{kt13} there should be no effect if the infinitesimal strain-free rigid body requirements are satisfied by the local coordinate system computations (for finite strains refer to Rankin \cite{Rankin}, where avoiding this assumption is accomplished).  Hence, adding $\mathbf{Q}^T\delta\mathbf{d}^L$ to \eqref{kt13} yields
\begin{equation}\label{kt14}
\delta\mathbf{d}^i_{\ell}=\mathbf{Q}^T\delta\mathbf{d}^{i}+\mathbf{Q}^T\left(\begin{array}{c}
y^{iL} \\ 
-x^{iL}
\end{array} \right)\delta\theta.
\end{equation}
To find $\delta\theta$,
take the derivative of \eqref{angle12} by using $\frac{d(\tan^{-1}u)}{dx}=\frac{1}{1+u^2}\frac{du}{dx}$.  This yields
\begin{equation}\label{kt15}
\delta\theta=\frac{1}{1+\frac{b^2}{a^2}}\delta(-ba^{-1})=\frac{a^2}{a^2+b^2}(-\delta ba^{-1}+a^{-2}b\delta a)
=\frac{a^2}{a^2+b^2}\left(\frac{b\delta a}{a^2}-\frac{a\delta b}{a^2}\right).
\end{equation}
Simplification and reorganizing \eqref{kt15} results in
\begin{equation}\label{kt16}
\delta\theta=\frac{b\delta a-a\delta b}{a^2+b^2}=\frac{1}{a^2+b^2}(b\mathbf{c}^T-a\mathbf{a}^T_{\ell})\delta \mathbf{d}
=\mathbf{v}^T\delta \mathbf{d}.
\end{equation}
Inserting $\delta \theta=\mathbf{v}^T\delta \mathbf{d}$ into \eqref{kt14} yields
\begin{equation}\label{kt17}
\delta\mathbf{d}^i_{\ell}=\mathbf{Q}^T\delta\mathbf{d}^{i}+\mathbf{Q}^T\left(\begin{array}{c}
y^{iL} \\ 
-x^{iL}
\end{array} \right)\mathbf{v}^T\delta \mathbf{d}.
\end{equation}
Now recognize that $\mathbf{Q}^T\left(\begin{array}{c}
y^{iL} \\ 
-x^{iL}
\end{array} \right)=\left(\begin{array}{c}
y^i_{\ell} \\ 
-x^i_{\ell}
\end{array} \right)$, and \eqref{kt17} becomes
\begin{equation}\label{kt18}
\delta\mathbf{d}^i_{\ell}=\mathbf{Q}^T\delta\mathbf{d}^{i}+\left(\begin{array}{c}
y^i_{\ell} \\ 
-x^i_{\ell}
\end{array} \right)\mathbf{v}^T\delta \mathbf{d}.
\end{equation}
Utilizing \eqref{kt18}, an alternative form is created for all element nodes as
\begin{equation}\label{kt19}
\delta \mathbf{d}_{\ell}=(\mathbf{\bar{Q}}+\hat{\mathbf{x}}_{\ell}\mathbf{v}^T)\delta \mathbf{d},
\end{equation}
where
\begin{equation*}\label{kt20}
\mathbf{\bar{Q}}=\left[\begin{array}{cccc}
[\mathbf{Q}^T] & \mathbf{0} & \hdots & \mathbf{0} \\ 
\mathbf{0} & [\mathbf{Q}^T] &  & \vdots\\ 
\vdots &  & \ddots & \vdots \\ 
\mathbf{0} & \hdots & \hdots & [\mathbf{Q}^T]
\end{array} \right], \quad \mathbf{0}=\left[\begin{array}{cc}
0 & \ \ 0 \\ 
0 & \ \ 0
\end{array} \right]
\end{equation*}
and
\begin{equation*}\label{kt21}
\hat{\mathbf{x}}^T_{\ell}=[\begin{array}{ccccccc}
y^1_{\ell} & -x^1_{\ell} & y^2_{\ell} & -x^2_{\ell} & \hdots & y^n_{\ell} & -x^n_{\ell}
\end{array}  ].
\end{equation*}
Observe that $\mathbf{\bar{Q}}$ is a $2n$ by $2n$ matrix.  Next, compare \eqref{kt19} with \eqref{kt3} to see that
\begin{equation}\label{kt22}
\mathbf{T}=\mathbf{\bar{Q}}+\hat{\mathbf{x}}_{\ell}\mathbf{v}^T.
\end{equation}

Finding the initial stiffness matrix $\mathbf{k}_{t\sigma}$ remains to complete construction of the tangent stiffness (see \eqref{kt6}).  See \eqref{kt5} and observe how the initial stiffness matrix arises from
\begin{equation}\label{kt23}
\delta\mathbf{T}^T\mathbf{q}_{\ell}=\mathbf{k}_{t\sigma}\delta\mathbf{d}.
\end{equation}
To find the first part of \eqref{kt23} take the variation of $\mathbf{T}^T$ to get
\begin{equation}\label{kt25}
\delta\mathbf{T}^T\mathbf{q}_{\ell}=\delta\mathbf{T}^1\mathbf{q}^1_{\ell}+\delta\mathbf{T}^2\mathbf{q}^2_{\ell}
+\hdots=\sum^{2n}_{j=1}\delta\mathbf{T}^j\mathbf{q}^j_{\ell},
\end{equation}
where $\mathbf{q}^j_{\ell}$ is the $j$th component of $\mathbf{q}_{\ell}$ (which is a scalar) and $\mathbf{T}^j$ is the $j$th column of $\mathbf{T}^T$.  Focus now on the first term of the summation \eqref{kt25} and take the transpose of \eqref{kt22} to obtain
\begin{equation}\label{kt26}
\delta\mathbf{T}^1\mathbf{q}^1_{\ell}=\mathbf{q}^1_{\ell}\delta\left\lbrace \left\lbrace  \begin{array}{c}
\mathbf{e}'_1 \\ 
\mathbf{0} \\ 
\vdots \\ 
\mathbf{0}
\end{array}\right\rbrace  +y^1_{\ell}\mathbf{v} \right\rbrace  =\mathbf{q}^1_{\ell}\mathbf{G}^1\delta\mathbf{d},
\end{equation}
where $\mathbf{0}^T=\begin{bmatrix}
0 & \ \ 0
\end{bmatrix}$.  Observe that $\mathbf{G}^1\delta\mathbf{d}$ is found from \eqref{kt26}.  This is written as
\begin{equation}\label{kt27}
\mathbf{G}^1\delta\mathbf{d}=\delta\left\lbrace \left\lbrace  \begin{array}{c}
\mathbf{e}'_1 \\ 
\mathbf{0} \\ 
\vdots \\ 
\mathbf{0}
\end{array}\right\rbrace  +y^1_{\ell}\mathbf{v} \right\rbrace  =\left\lbrace  \begin{array}{c}
\mathbf{e}'_2 \\ 
\mathbf{0} \\ 
\vdots \\ 
\mathbf{0}
\end{array}\right\rbrace\delta\theta+\delta y^1_{\ell}\mathbf{v}+y^1_{\ell}\delta \mathbf{v}.
\end{equation}
Next, realize that \eqref{kt19} yields $\delta y^1_{\ell}$, i.e.,
\begin{equation}\label{kt28}
\delta y^1_{\ell}=\left \lbrace \left[\begin{array}{ccccc}
\mathbf{e}'^T_2 &\mathbf{0}  & \mathbf{0} & \hdots & \mathbf{0}
\end{array} \right] - x^1_{\ell}\mathbf{v}^T\right\rbrace\delta\mathbf{d}.
\end{equation}
Recognize this by considering the generic variable $\mathbf{w}$.  Then realize that the variation of $\mathbf{w}$ in local coordinates is expressible in terms of $\mathbf{w}$ in global coordinates as (see \eqref{kt19})
\begin{equation}\label{kt29}
\delta \mathbf{w}_{\ell}=(\mathbf{\bar{Q}}+\hat{\mathbf{x}}_{\ell}\mathbf{v}^T)\delta \mathbf{w},
\end{equation}
where $\mathbf{w}=\mathbf{X}+\mathbf{d}$ and $\mathbf{w}_{\ell}=\mathbf{X}_{\ell}+\mathbf{d}_{\ell}=\mathbf{x}_{\ell}$.  Specifically, $\mathbf{w}^T_{\ell}=\left \lbrace\begin{array}{ccccc}
 x^1_{\ell}&y^1_{\ell}  &\hdots  &x^n_{\ell}  &y^n_{\ell}\end{array} \right\rbrace $.  Then recognize that $\delta\mathbf{w}=\delta\mathbf{d}$ because $\delta\mathbf{X}=0$.  Equation \eqref{kt28} is found from the row of \eqref{kt29} corresponding to the variation $\delta y^1_{\ell}$.

For now ignore the last term of \eqref{kt27}, and use \eqref{kt16}, \eqref{kt27} and \eqref{kt28} to get
\begin{equation}\label{kt30}
\mathbf{G}^{1,a}=\left \lbrace\begin{array}{c}
\mathbf{e}'_2 \\ 
\mathbf{0} \\ 
\vdots \\ 
\mathbf{0}
\end{array} \right \rbrace \mathbf{v}^T+\mathbf{v}\left \lbrace \begin{array}{c}
\mathbf{e}'_2 \\ 
\mathbf{0} \\ 
\vdots \\ 
\mathbf{0}
\end{array}\right \rbrace ^T-x^1_{\ell}\mathbf{v}\mathbf{v}^T,
\end{equation}
which is symmetric.  The variation $\delta\mathbf{v}$ is required to obtain the complete form of $\mathbf{G}^1$.  With this in mind, use \eqref{kt16} to write 
\begin{equation}\label{kt31}
\mathbf{v}=\frac{1}{a^2+b^2}(b\mathbf{c}-a\mathbf{a}_{\ell}).
\end{equation}
Then, taking the variation of \eqref{kt31} and using algebra \cite{yaw}, eventually, one finds
\begin{equation}\label{kt35}
\delta\mathbf{v} = 
\left[\frac{2ab(\mathbf{a_{\ell}a_{\ell}}^T-\mathbf{cc}^T)+(a^2-b^2)(\mathbf{ca_{\ell}}^T+\mathbf{a_{\ell}c}^T)}{(a^2+b^2)^2}\right]\delta\mathbf{d}
\equiv \mathbf{V}^T\delta\mathbf{d},
\end{equation}
where $\vm{V}$ is a symmetric matrix.  The 
matrix $\mathbf{V}^T$ is as determined by Crisfield and Moita~\cite{Crisfield3} yet their denominator is not squared (a likely typographical error).  Note also, according to a study by \cite{Crisfieldv2,Rankin}, $\mathbf{G}^{1,b}$ below, which contains the variation of $\mathbf{v}$, can be discarded, as its affect on convergence is insignificant.  Nevertheless, it is kept here to be complete.  Therefore, with
\begin{equation}\label{kt36}
\mathbf{G}^{1,b}=y^1_{\ell}\mathbf{V}^T,
\end{equation}
the full version of $\mathbf{G}^1$ is finally found as
\begin{equation}\label{kt37}
\mathbf{G}^1=\mathbf{G}^{1,a}+\mathbf{G}^{1,b}.
\end{equation}
However, recall that $\mathbf{G}^1$ only provides the first term in summation \eqref{kt25}.  Other matrices $\mathbf{G}^j$ are similarly found.  As a result, $\mathbf{k}_{t\sigma}=\sum^{2n}_{j=1}\mathbf{q}^j_{\ell}\mathbf{G}^j$ provides the expression for the initial stiffness and then the total tangent stiffness matrix is given according to~\eqref{kt6}.

Additionally, $\mathbf{G}^2$ is 
provided.  The second term in the summation of \eqref{kt25} results in
\begin{equation}\label{kt38}
\delta\mathbf{T}^2\mathbf{q}^2_{\ell}=\mathbf{q}^2_{\ell}\delta\left\lbrace \left\lbrace  \begin{array}{c}
\mathbf{e}'_2 \\ 
\mathbf{0} \\ 
\vdots \\ 
\mathbf{0}
\end{array}\right\rbrace  -x^1_{\ell}\mathbf{v} \right\rbrace  =\mathbf{q}^2_{\ell}\mathbf{G}^2\delta\mathbf{d}.
\end{equation}
From \eqref{kt38}, $\mathbf{G}^2\delta\mathbf{d}$ is found, which is written as
\begin{equation}\label{kt39}
\mathbf{G}^2\delta\mathbf{d}=\delta\left\lbrace \left\lbrace  \begin{array}{c}
\mathbf{e}'_2 \\ 
\mathbf{0} \\ 
\vdots \\ 
\mathbf{0}
\end{array}\right\rbrace  -x^1_{\ell}\mathbf{v} \right\rbrace=\left\lbrace  \begin{array}{c}
-\mathbf{e}'_1 \\ 
\mathbf{0} \\ 
\vdots \\ 
\mathbf{0}
\end{array}\right\rbrace\delta\theta+\delta(-x^1_{\ell})\mathbf{v}+(-x^1_{\ell})\delta \mathbf{v}.
\end{equation}
The variation $\delta x^1_{\ell}$ arises in that same way as shown after \eqref{kt28}, i.e.,
\begin{equation}\label{kt40}
\delta x^1_{\ell}=\left \lbrace \left[\begin{array}{ccccc}
\mathbf{e}'^T_1 &\mathbf{0}  & \mathbf{0} & \hdots & \mathbf{0}
\end{array} \right] +y^1_{\ell}\mathbf{v}^T\right\rbrace\delta\mathbf{d}.
\end{equation}
Equation \eqref{kt39} and use of \eqref{kt16}, \eqref{kt40} and \eqref{kt35} yields
\begin{equation*}\label{kt41}
\mathbf{G}^2=\left \lbrace\begin{array}{c}
-\mathbf{e}'_1 \\ 
\mathbf{0} \\ 
\vdots \\ 
\mathbf{0}
\end{array} \right \rbrace \mathbf{v}^T+\mathbf{v}\left \lbrace \begin{array}{c}
-\mathbf{e}'_1 \\ 
\mathbf{0} \\ 
\vdots \\ 
\mathbf{0}
\end{array}\right \rbrace ^T-y^1_{\ell}\mathbf{v}\mathbf{v}^T-x^1_{\ell}\mathbf{V}^T.
\end{equation*}
Finally, the generic cases of $\mathbf{G}^{2i-1}$ and $\mathbf{G}^{2i}$ are provided.  In general, for~$i=1$~to~$n$
\begin{equation*}\label{kt42}
\mathbf{G}^{2i-1}=\begin{array}{c}
1 \\ 
\vdots \\ 
i \\ 
\vdots \\ 
n
\end{array} \left \lbrace\begin{array}{c}
\mathbf{0} \\
\vdots\\ 
\mathbf{e}'_2 \\ 
\vdots \\ 
\mathbf{0}
\end{array} \right \rbrace \mathbf{v}^T+\mathbf{v}\left \lbrace \begin{array}{c}
\mathbf{0} \\ 
\vdots \\
\mathbf{e}'_2 \\ 
\vdots \\ 
\mathbf{0}
\end{array}\right \rbrace ^T-x^i_{\ell}\mathbf{v}\mathbf{v}^T+y^i_{\ell}\mathbf{V}^T,
\end{equation*}
\begin{equation*}\label{kt43}
\mathbf{G}^{2i}=\left \lbrace\begin{array}{c}
\mathbf{0} \\
\vdots\\ 
-\mathbf{e}'_1 \\ 
\vdots \\ 
\mathbf{0}
\end{array} \right \rbrace \mathbf{v}^T+\mathbf{v}\left \lbrace \begin{array}{c}
\mathbf{0} \\ 
\vdots \\
-\mathbf{e}'_1 \\ 
\vdots \\ 
\mathbf{0}
\end{array}\right \rbrace ^T-y^i_{\ell}\mathbf{v}\mathbf{v}^T-x^i_{\ell}\mathbf{V}^T.
\end{equation*}

\subsection{Element Internal Forces}
Within a (nonlinear) co-rotational analysis, element local internal forces are needed for an implicit analysis with iterations for equilibrium.  For linear elastic problems it is sufficient to use, for a single element~$E$,
\begin{equation}\label{eif1}
	\mathbf{q}^E_{\ell int}=\mathbf{k}_{t\ell}\mathbf{d}_{\ell}=\mathbf{k}_{E}\mathbf{d}_{\ell}.   \quad \quad [2n_v \times 1]
\end{equation}
For plastic problems it is better to directly use the results from the stress integration process and also correct for the adjustment provided by the stabilization matrix.  For a single element~$E$,
\begin{equation}\label{eif1b}
	\mathbf{q}^E_{\ell int}=A_E \ t \ \mathbf{B}^T\bsym{\sigma}+\mathbf{k}_E^s \mathbf{d}_{\ell},
\end{equation}
where $A_E$ is the area of polygon element $E$ and $\bsym{\sigma}=[\sigma_x \ \sigma_y \ \sigma_{xy}]^T$ is the vector of engineering stresses, constant across the given element for a particular time step in the nonlinear analysis.
Then similar to FEM the individual transformed internal force vectors for all elements are assembled into the global internal force vector using the assembly operator \cite{Hughes}.  That is,
\begin{equation}\label{eif2}
	\mathbf{F}_{int}=\overset{n_{elem}}{\underset{E=1}{\mathbf{\mathsf{A}}}}\mathbf{T}^T\mathbf{q}^{E}_{\ell int}.
\end{equation}

\noindent
\textbf{Remarks}
\begin{enumerate}[(i)]
	\item The necessary pieces for an implicit nonlinear analysis are now in place.
	\item The tangent stiffness for an individual element is represented by $\mathbf{k}_T$ in \eqref{kt6}.  Using this information a global structure tangent stiffness matrix, $\mathbf{K}_T$, is assembled, recognizing that the number of element degrees of freedom vary depending on the number of vertices in a given polygon element.
	\item The global vector of internal forces is prescribed by \eqref{eif2}.
	\item The global vector of external forces are prescribed by the analyst and constructed in the standard way like FEM per the linear form \eqref{ext1}.
	\item It is now possible to calculate a residual, $\mathbf{g}=\mathbf{F}_{int}-\lambda \mathbf{F}_{ext}$, for a given load or iteration step.  Here, $\lambda$ is the arc length control loading parameter, and $\mathbf{g}$ is the residual vector during iterations for equilibrium.
\end{enumerate}

\section{Co-rotational VEM - Plasticity}\label{sec:corotplast}
A co-rotational analysis is a nonlinear analysis that includes geometric nonlinearities.  Additionally, material nonlinearities such as 2D plane stress plasticity are easily included.  For such a case material properties must be updated during each load step.  The local VEM stiffness is the same as the case for elasticity except that a consistent elasto-plastic modular matrix, $\mathbf{C}_{ep}$, is inserted into the expressions for $\mathbf{k}_E^c$ and $\mathbf{k}_E^s$.  Then the matrix $\mathbf{k}_{t\ell}$ becomes
\begin{equation}
	\mathbf{k}_{t\ell}=\mathbf{k}_E= \mathbf{k}_E^c+\mathbf{k}_E^s=tA_E\mathbf{B}^T \mathbf{C}_{ep} \mathbf{B}+\mathbf{k}_E^s(C_{ep}).
\end{equation}
As strains evolve during each load step, stresses and $\mathbf{C}_{ep}$ are updated according to the $J2$ plasticity formulation with radial return (see Simo and Taylor \cite{Simo1} and Simo and Hughes \cite{Simo2}).  All other formulas for VEM and co-rotation remain the same.

\section{Numerical Implementation}\label{sec:numimplementation}
For the nonlinear analysis procedure implicit Newton-Raphson iterations are used to enforce global equilibrium.  In addition, for problems involving plasticity, implicit Newton-Raphson iterations are used with radial return at the constitutive level for each element.  Global equilibrium is enforced by an arc-length path following scheme \cite{Crisfieldv1}.  As a result, the nonlinear co-rotational analysis is carried out with the ingredients described in prior sections: (i) a path following scheme (arc-length method), (ii) global external force vector, (iii) global internal force vector, and (iv) consistent global tangent stiffness matrix.

As alluded to previously, in section \ref{sec:angle}, the co-rotational formulation requires shape function derivatives at the centroid of each element.  Initially, it was thought that these terms needed to be calculated by some `alternative' means.  Derivatives of moving least squares (MLS) \cite{tb:meshless} shape functions and mean value coordinates (MVC) \cite{Floater} shape functions were both explored and shown to work well.  Yet, it is preferrable to avoid these alternatives and recognize that the approximations of VEM shape function derivatives are available in the strain displacement matrix,
\begin{equation}
		\mathbf{B}=\left[ \begin{array}{ccccccc}
		\partial_{x}\phi_1 & 0 & \partial_{x}\phi_2 & 0 & ... & \partial_{x}\phi_{n_v} & 0 \\
		0 & \partial_{y}\phi_1 & 0 & \partial_{y}\phi_2 &... & 0 & \partial_{y}\phi_{n_v} \\
		\partial_{y}\phi_1 &\partial_{x}\phi_1 & \partial_{y}\phi_2 & \partial_{x}\phi_2 & ... & \partial_{y}\phi_{n_v} &\partial_{x}\phi_{n_v} 
		\end{array}\right].
\end{equation}  
Although these derivatives (needed at element centroid) are constant across each element, numerical results verify that they are as effective as the other methods investigated.  Furthermore, extracting and using these derivatives keeps the formulation within the confines of the virtual element method.

\section{Numerical Results}\label{sec:numresults}
A variety of simulations are provided to illustrate the utility of the co-rotational VEM.  In particular, representative static elastic and plastic simulations with results are provided.  All simulations are accomplished by way of implicit nonlinear analysis. Newton-Raphson iterations are used to enforce equilibrium using arc length control\cite{Crisfieldv1}.  All polygonal mesh generation is accomplished by using Polymesher \cite{talischi}.

\subsection{Linear Elastic Cantilever Beam}
\begin{figure}
\centering
\mbox{
\subfigure[]{\epsfig{file=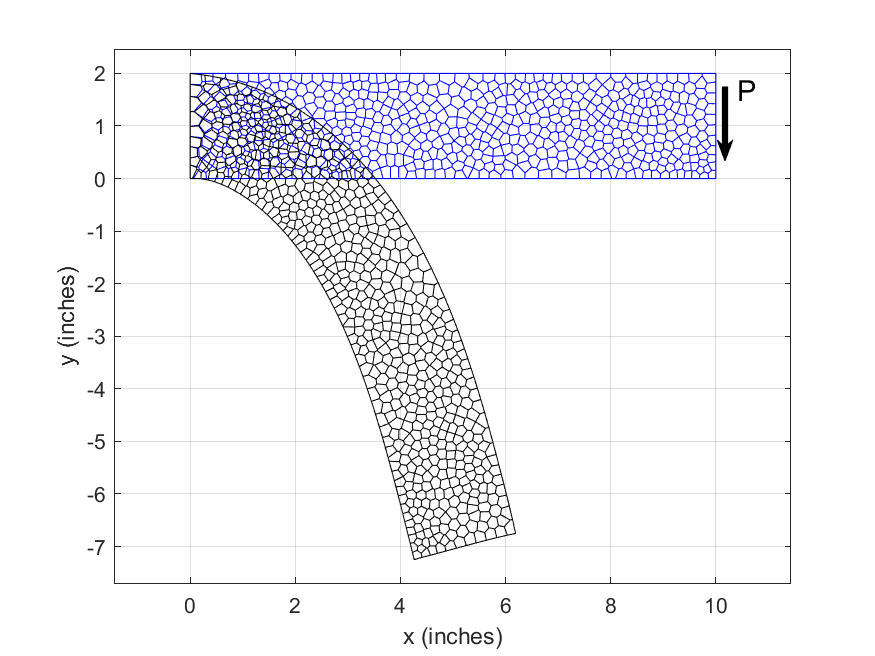,width=0.48\textwidth}}
\subfigure[]{\epsfig{file=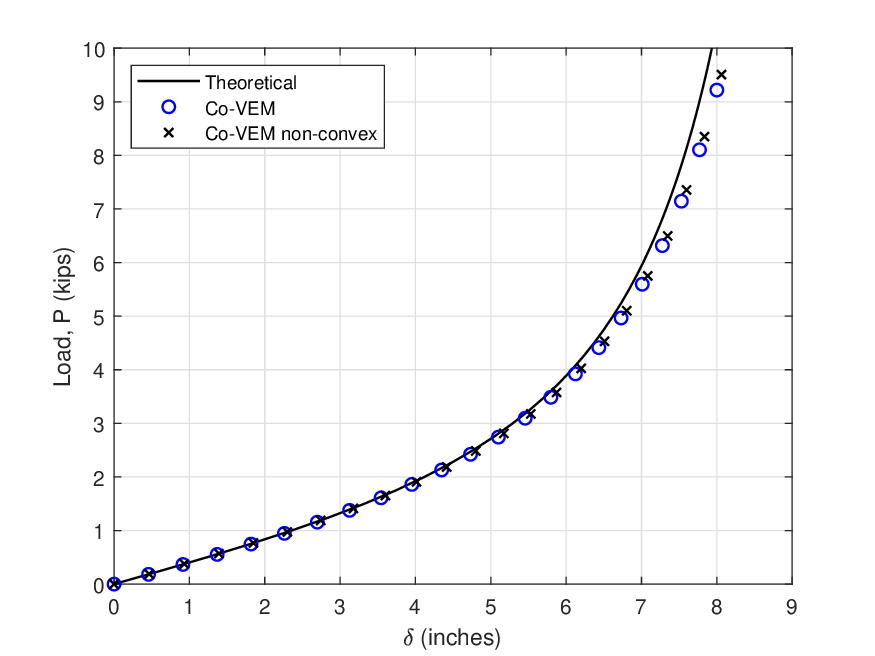,width=0.48\textwidth}}
}

\mbox{
\subfigure[]{\epsfig{file=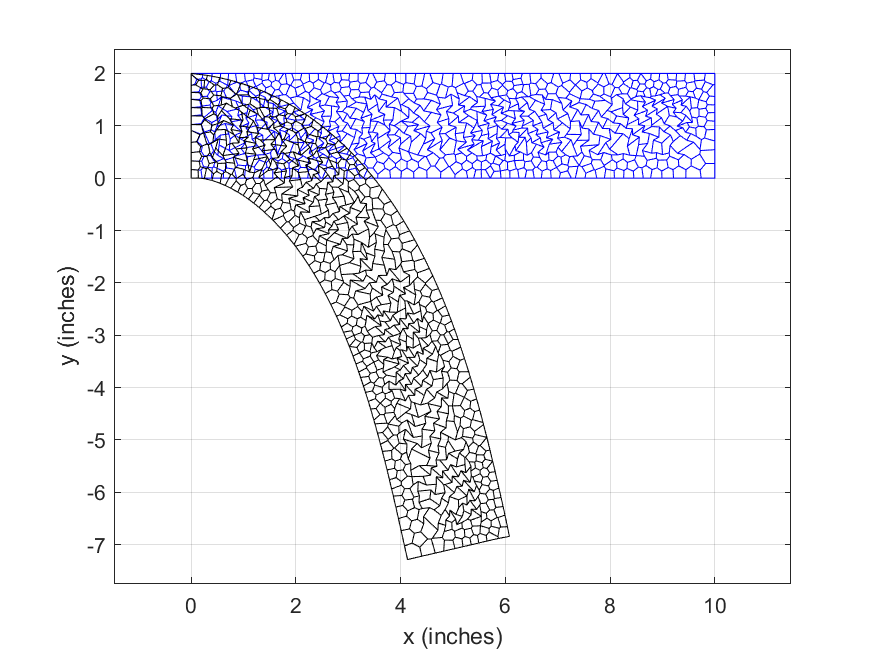,width=0.48\textwidth}}
\subfigure[]{\epsfig{file=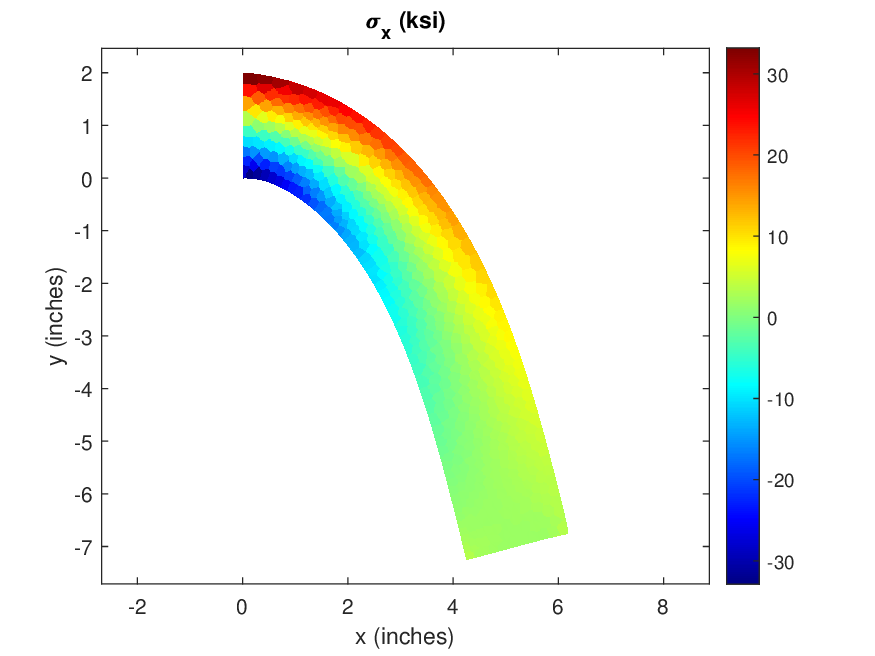,width=0.48\textwidth}}
}

\mbox{
\subfigure[]{\epsfig{file=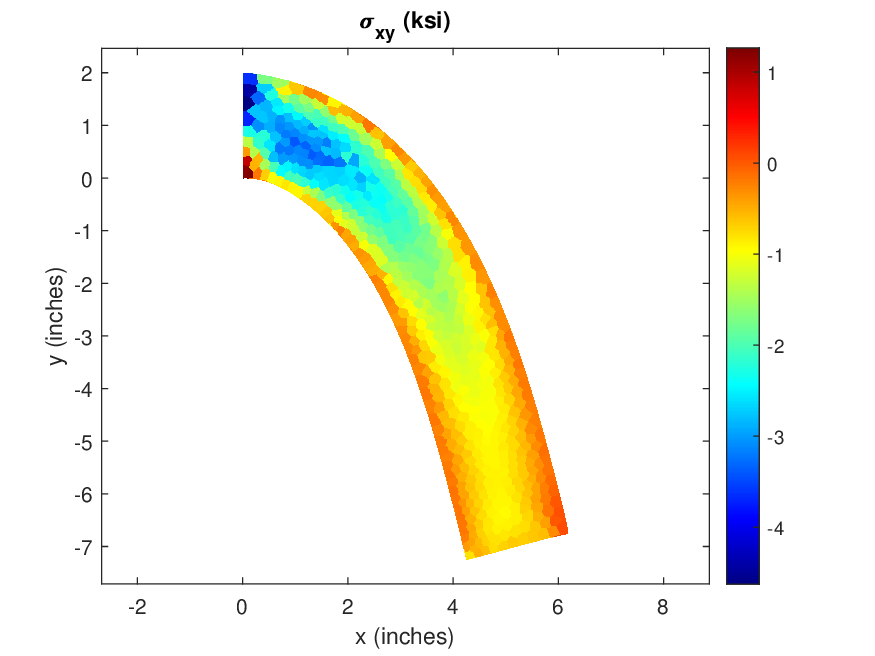,width=0.48\textwidth}}
}

\vspace*{-0.1in}
\caption{2D Linear Elastic Cantilever:  (a)~VEM undeformed and deformed structure (620 elements); (b)~Load versus displacement results compared to theoretical solution; (c)~VEM results containing non-convex elements ; (d)~$\sigma_x$ stress; and (e)~$\sigma_{xy}$ stress.}\label{numericalcant}
\end{figure}

For the condition of plane stress a linear elastic cantilever is loaded downward at its free end.  The cantilever has thickness $t=2$~inches, modulus of elasticity $E_Y=100$~ksi and Poisson's ratio $\nu=0.0$.  The theoretical solution, for load versus displacement, including bending and axial deformations is found in \cite{yaw}.  In Figure~\ref{numericalcant}b, for models with convex or non-convex polygons, load versus displacement results are shown to be in excellent agreement with the theoretical solution.  At higher loads the slight discrepancy is due to shear deformations being absent from the theoretical solution.  The cantilevers modeled in this example and shown in Figures~\ref{numericalcant}a and~c are modeled with 620 polygons.  The stresses shown are plotted in local co-rotated coordinates for each polygon cell of the model.

\subsection{Linear Elastic Ring}
\begin{figure}
\centering
\mbox{
\subfigure[]{\epsfig{file=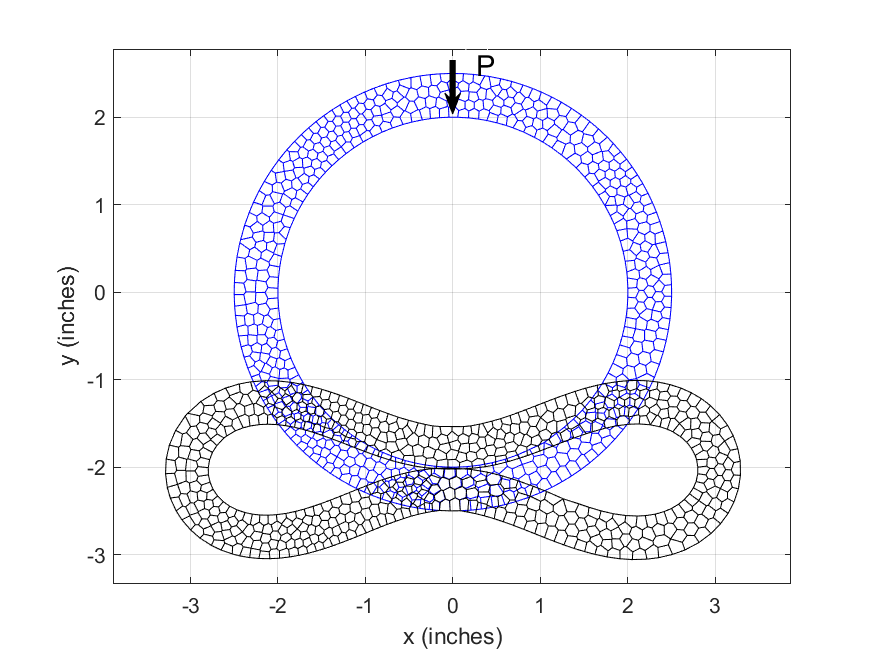,width=0.48\textwidth}}
\subfigure[]{\epsfig{file=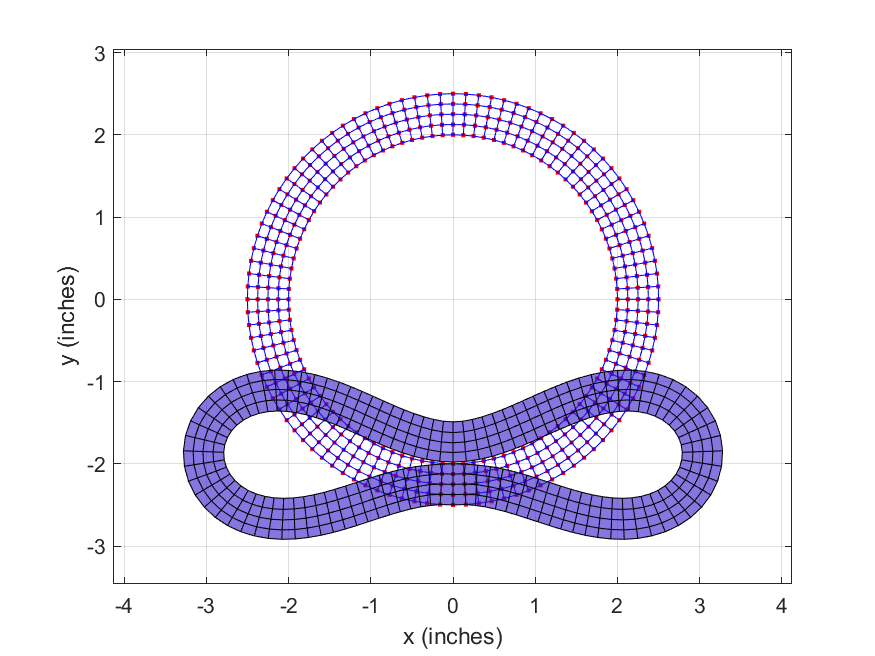,width=0.48\textwidth}}
}

\mbox{
\subfigure[]{\epsfig{file=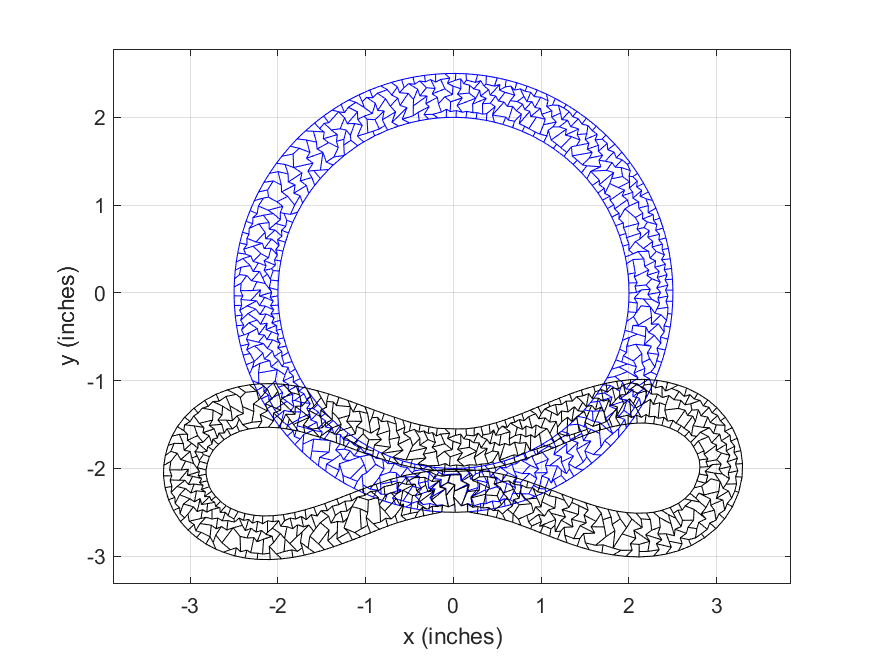,width=0.48\textwidth}}
\subfigure[]{\epsfig{file=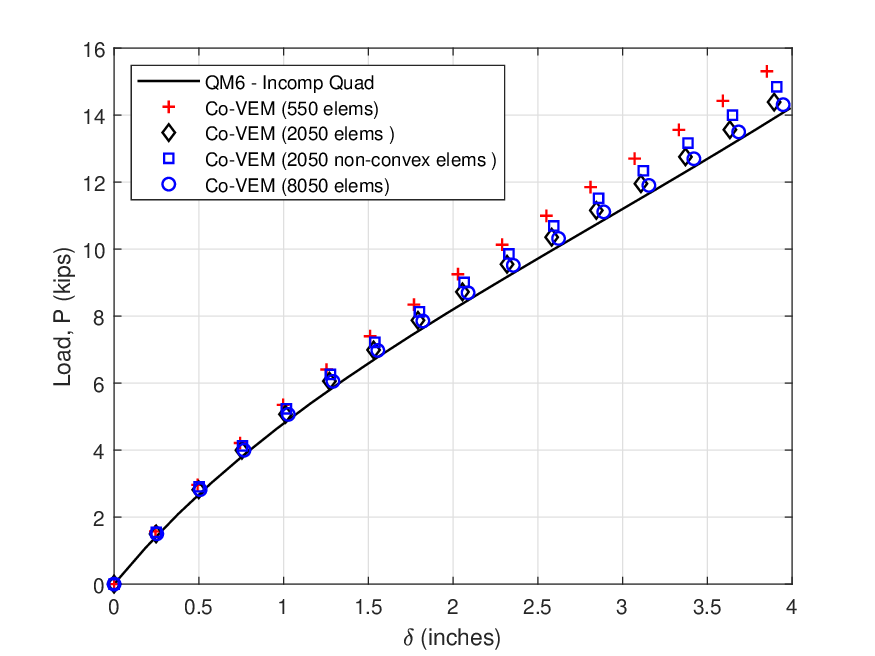,width=0.48\textwidth}}
}

\mbox{
\subfigure[]{\epsfig{file=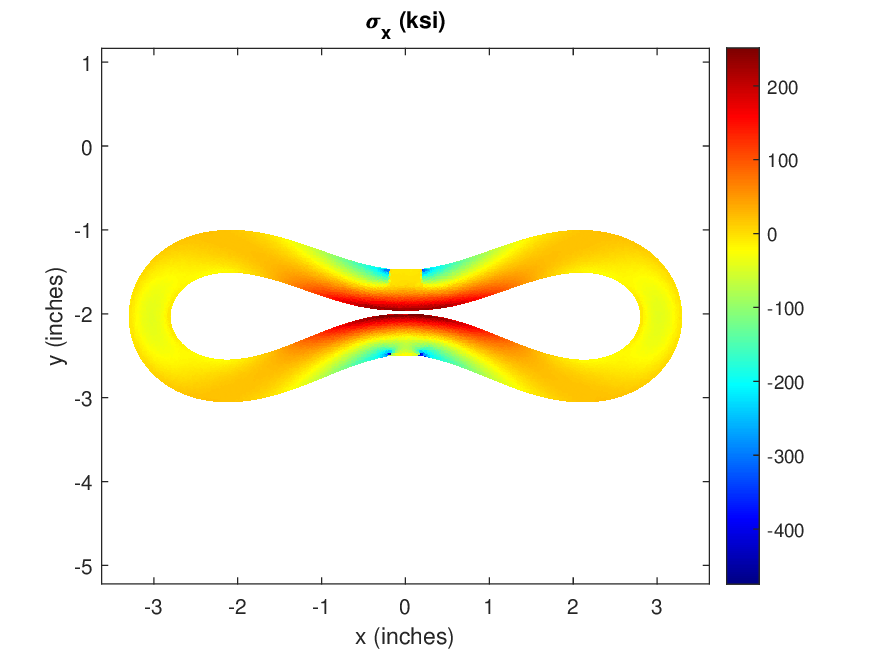,width=0.48\textwidth}}
\subfigure[]{\epsfig{file=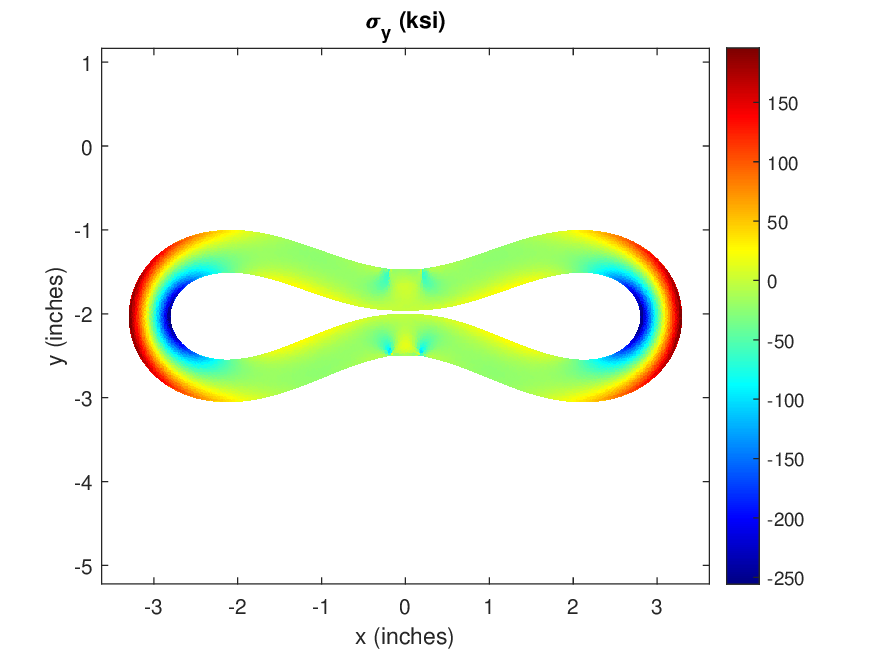,width=0.48\textwidth}}
}

\vspace*{-0.1in}
\caption{2D Linear Elastic Ring:  (a)~VEM undeformed and deformed structure (550 elements); (b)~FEM incompatible modes element results (400 elements); (c)~VEM results with 550 non-convex elements ; (d)~Load versus displacement convergence study;
 (e)~$\sigma_x$ stress; and (f)~$\sigma_y$ stress.}\label{numericalring}
\end{figure}

For the condition of plane stress a linear elastic ring is supported at the bottom and loaded at the top of the ring.  The inner and outer radii of the ring are $r_i=2.0$ and $r_o=2.5$~inches, respectively.  The ring has thickness $t=1$~inch, modulus of elasticity $E_Y=1000$~ksi and Poisson's ratio $\nu=0.3$.  In Figure~\ref{numericalring}d, refinement of the VEM mesh converges to similar results for 400 FEM co-rotational incompatible modes quadrilateral (QM6) elements \cite{WilsonQM6,TaylorQM6,Crisfield3}.  The FEM results for 816 elements have negligible difference verifying that 400 elements are adequate for a benchmark solution.  Significantly fewer QM6 elements are required since the incompatible modes elements are constructed to include bending mode enhancement (quadratic polynomials), whereas the VEM elements used herein are only based on first order polynomials.  Nevertheless, the co-rotational VEM formulation achieves similar results with appropriate mesh refinement.  Non-convex elements are confirmed to provide similar results in both Figures~\ref{numericalring}c and~d.

\subsection{Linear Elastic Arch} 
\begin{figure}
\centering
\mbox{
\subfigure[]{\epsfig{file=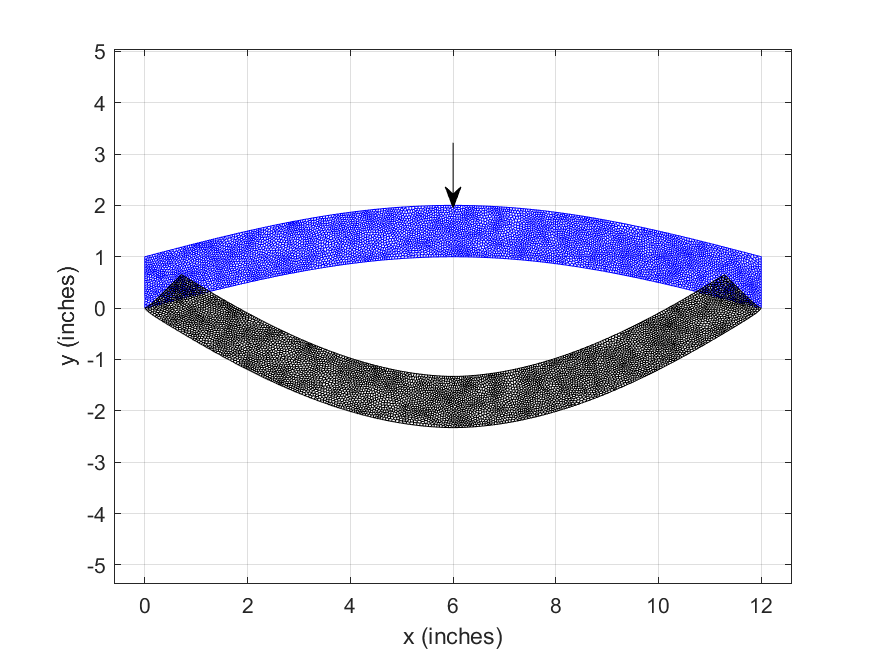,width=0.48\textwidth}}
\subfigure[]{\epsfig{file=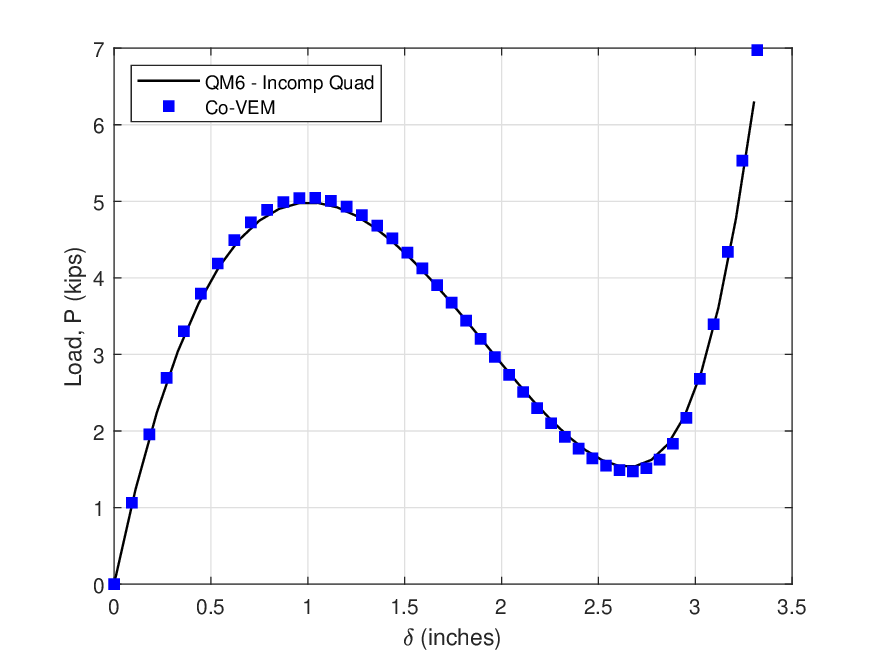,width=0.48\textwidth}}
}

\mbox{
\subfigure[]{\epsfig{file=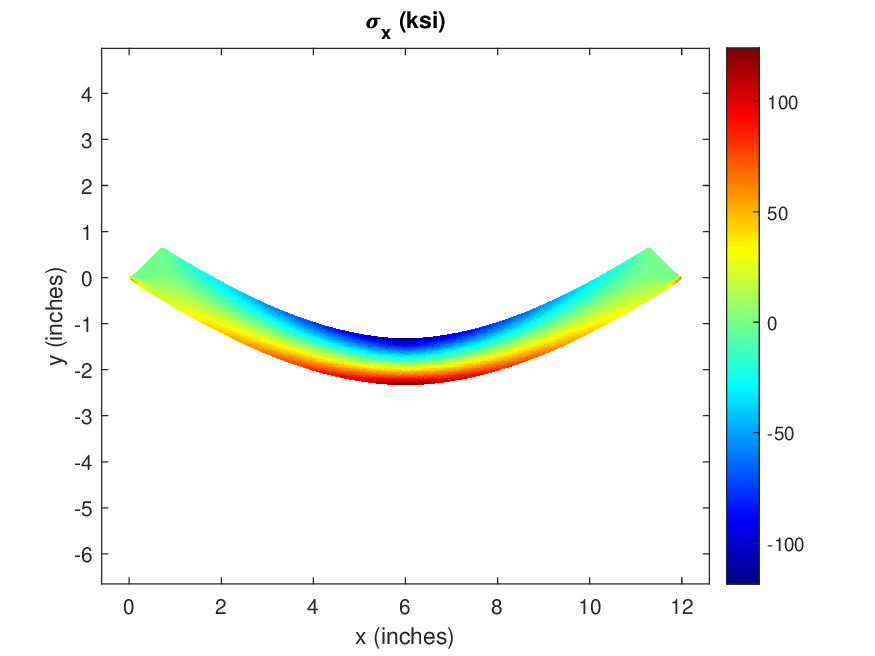,width=0.48\textwidth}}
}

\mbox{
}

\vspace*{-0.1in}
\caption{2D Linear Elastic Arch:  (a)~VEM undeformed and deformed structure (6000 elements); (b)~Load versus displacement compared to benchmark;
 (c)~$\sigma_x$ stress.}\label{numericalarch}
\end{figure}

In Figure~\ref{numericalarch}, a plane stress linear elastic arch, modeled with 6000 polygonal elements, is loaded at midspan.  The arch has thickness $t=1$~inch, span $L=12$~inches, modulus of elasticity $E_Y=1000$~ksi, and Poisson's ratio $\nu=0.3$.  The geometry of the arch is constructed by first creating a rectangular discretization of $x=12$ by $y=1$, then for each y coordinate setting $y=y+\sin{\frac{\pi x}{L}}$.  The load displacement results show very good agreement with results for a FEM co-rotational incompatible modes quadrilateral (QM6) \cite{WilsonQM6,TaylorQM6,Crisfield3} arch composed of 1260 elements.  Pinned supports are provided at the left and right bottom corners of the arch.

\subsection{Plastic Cantilever Beam} 
\begin{figure}
\centering
\mbox{
\subfigure[]{\epsfig{file=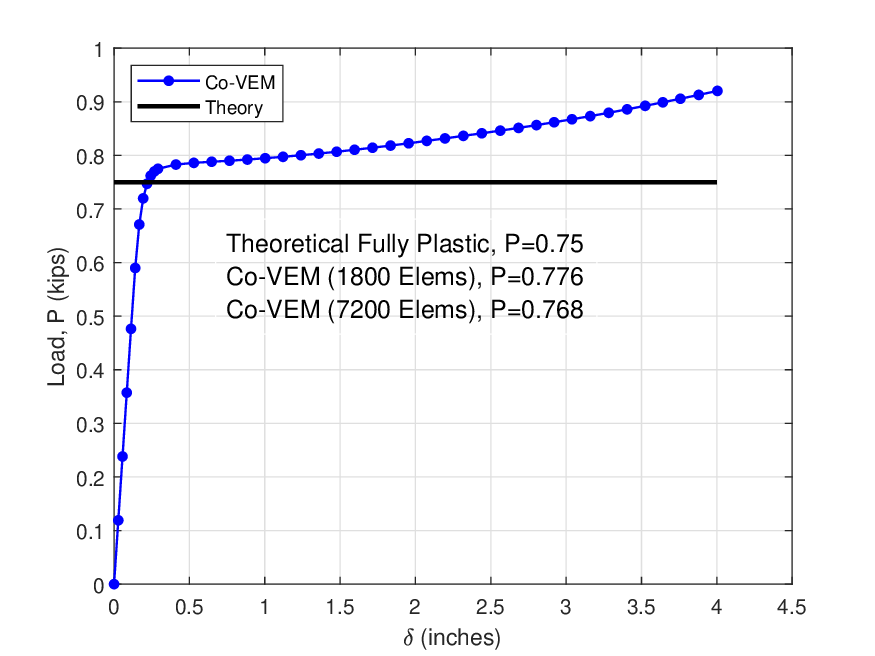,width=0.48\textwidth}}
\subfigure[]{\epsfig{file=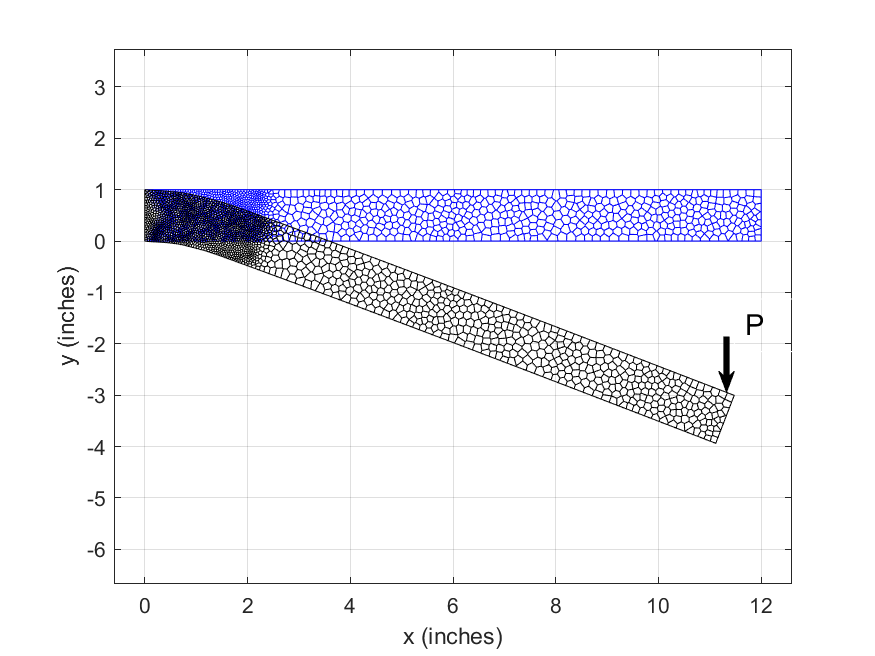,width=0.48\textwidth}}
}

\mbox{
\subfigure[]{\epsfig{file=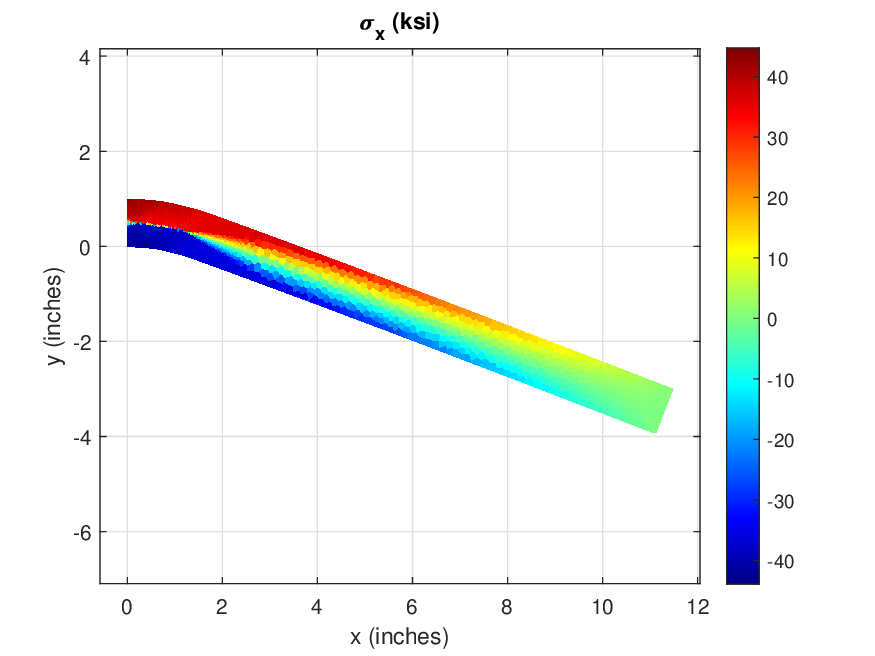,width=0.48\textwidth}}
\subfigure[]{\epsfig{file=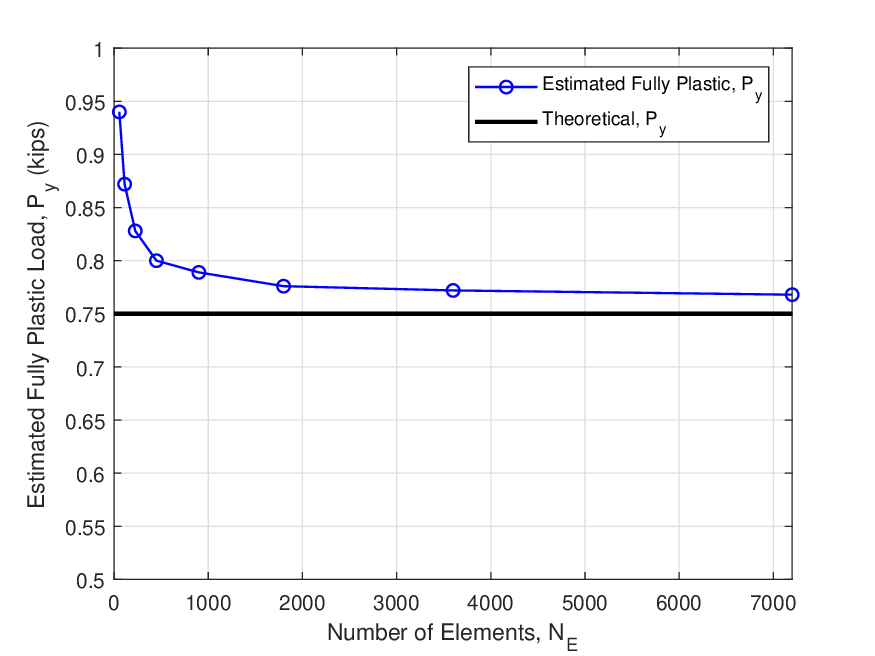,width=0.48\textwidth}}
}


\vspace*{-0.1in}
\caption{2D Elasto-Plastic Cantilever:  (a)~Load versus displacement; (b)~VEM undeformed and deformed configuration (1800 elements); (c)~$\sigma_x$ stress; and (d)~convergence study for fully plastic cross-section of cantilever compared to theoretical value.}\label{numericalplasticcant}
\end{figure}

For the condition of plane stress an elasto-plastic cantilever with 1800 elements is loaded downward at its free end (Figure~\ref{numericalplasticcant}).  The density of the elements is biased toward the cantilever's fixed support, as shown in Figure~\ref{numericalplasticcant}b, to capture the anticipated concentration of plastic flow near the support.  As expected a plastic hinge develops in the cantilever near the support (Figure~\ref{numericalplasticcant}c).  The cantilever has thickness $t=1$~inch, length $L=12$ inches, modulus of elasticity $E_Y=29000$~ksi, Poisson's ratio $\nu=0.3$, yield stress $\sigma_{yield}=36$~ksi, and linear hardening modulus $E_h=1$ ksi.  The theoretical load to cause a fully plastic hinge through the cantilever cross-section is $P_y=0.75$~kips.  With mesh refinement the numerically estimated load for fully plastic cross-section approaches the theoretical value as shown in Figure~\ref{numericalplasticcant}d.  The numerically estimated value of $P_y$ is 0.768~kips for a model with 7200 polygon elements.  For higher displacements the load displacement curve begins to climb due to tension stiffening as the cantilever rotates.  This expected behavior is captured by the co-rotational formulation and revealed in the load versus displacement plot of Figure~\ref{numericalplasticcant}a.

\subsection{Plastic Ring}\label{subsec:explicit_dynamics}
\begin{figure}
\centering
\mbox{
\subfigure[]{\epsfig{file=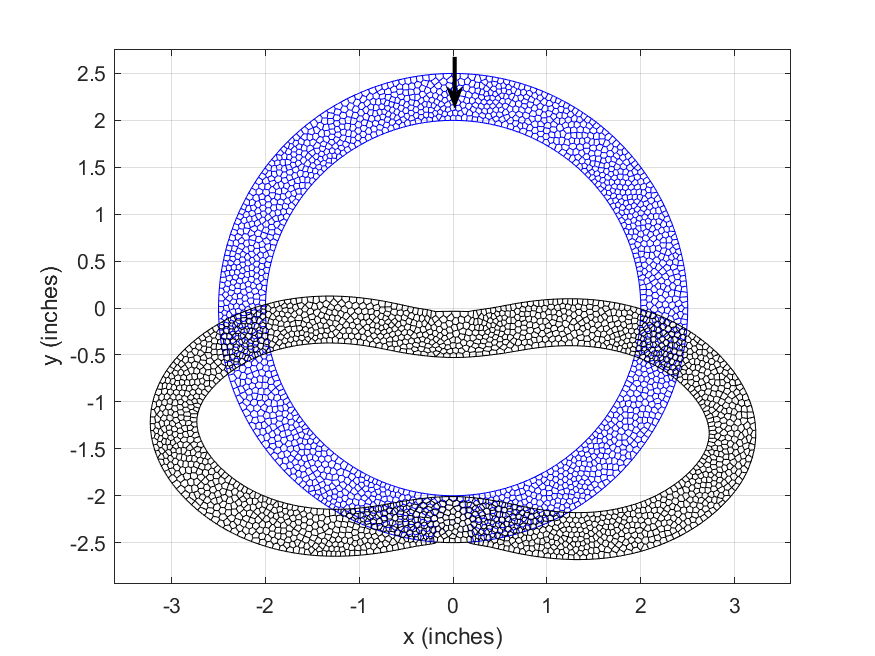,width=0.48\textwidth}}
\subfigure[]{\epsfig{file=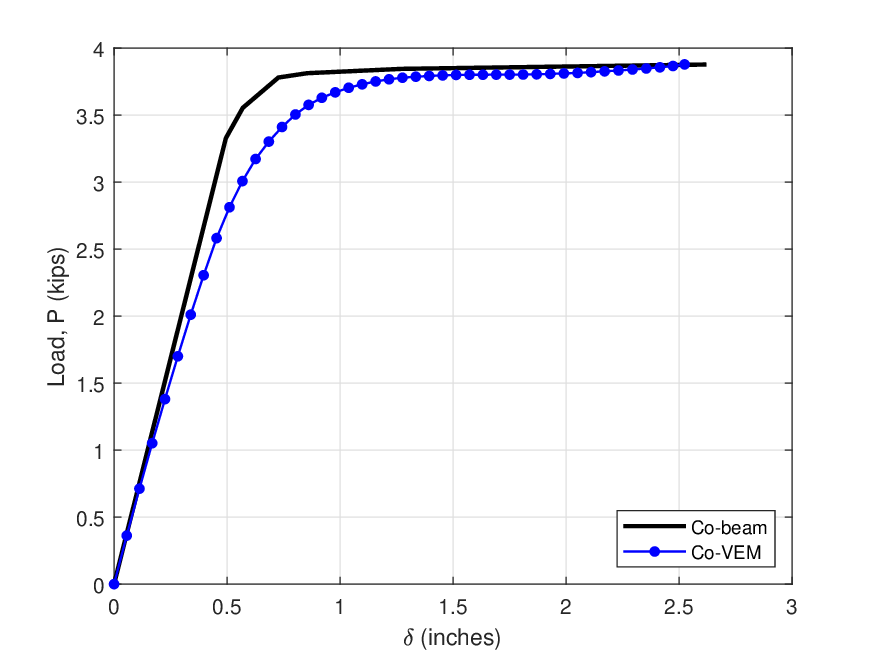,width=0.48\textwidth}}
}

\mbox{
\subfigure[]{\epsfig{file=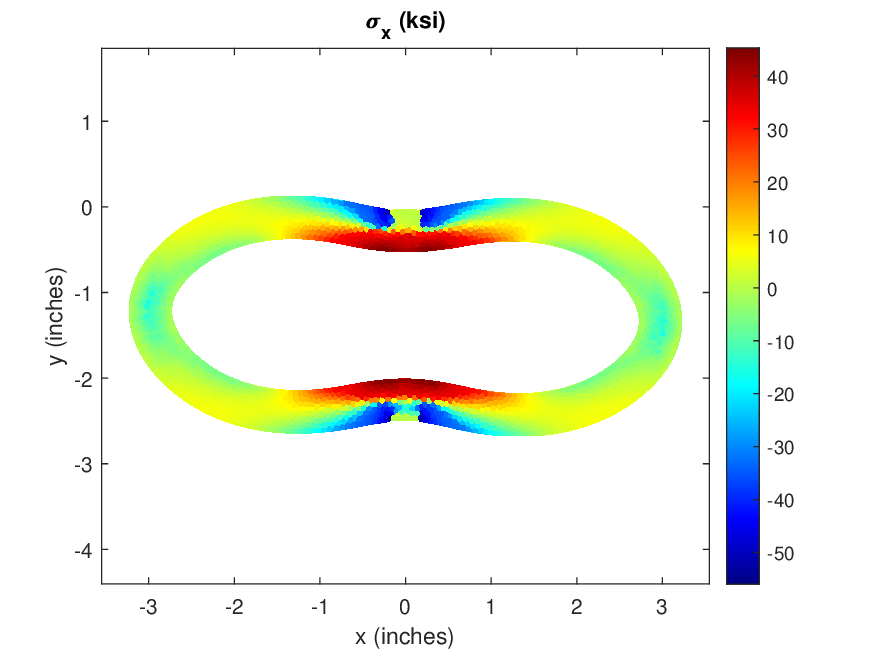,width=0.48\textwidth}}
\subfigure[]{\epsfig{file=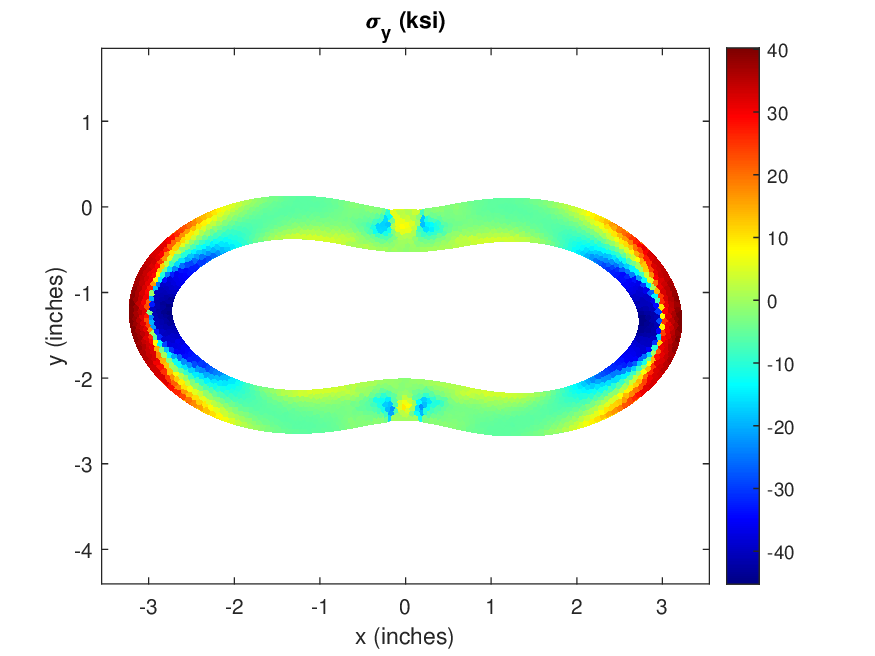,width=0.48\textwidth}}
}


\vspace*{-0.1in}
\caption{2D Elasto-Plastic Ring:  (a)~VEM undeformed and deformed structure (2050 elements); (b)~Load versus displacement for the ring structure; (c)~$\sigma_x$ stress; and (d)~$\sigma_y$ stress.}\label{numericalplasticring}
\end{figure}

For the case of plane stress an elasto-plastic ring with 2050 elements is loaded downward at the top in the middle (Figure~\ref{numericalplasticring}a).  The ring is supported in the middle at the base by a narrow region of nodes.  Clear plastic yielding is revealed by the nonlinear behavior of load versus displacement in Figure~\ref{numericalplasticring}b and the plastic hinges evident in Figures~\ref{numericalplasticring}c and d.  The load displacement results are in good agreement with a co-rotational elasto-plastic beam formulation.  Stresses at the final load and displacement are shown in Figures~\ref{numericalplasticring}c and~d.  The ring has thickness $t=1$~inch, modulus of elasticity $E_Y=29000$~ksi, Poisson's ratio $\nu=0.3$, yield stress $\sigma_{yield}=36$~ksi, and linear hardening modulus $E_h=1$~ksi.  The inner and outer radii of the ring are $r_i=2.0$ and $r_o=2.5$~inches, respectively.

\section{Conclusions}\label{sec:conclusions}
In this work, a first order virtual element method is included within a nonlinear co-rotational formulation.  The formulation, so constructed, is limited to small strains for elastic or plastic problems, but allows for large displacements and large rotations.  Capabilities of the formulation are illustrated by solving example problems and comparing them to theoretical results or the results obtained by finite element solutions.  It is evident from the results that co-rotational VEM for elastic and plastic problems is a viable scheme for solving nonlinear problems.

All numerical methods have limitations.  The formulation herein is no different and has typical limitations observed as follows:  (a) element instability beyond small to moderate strains, (b) need for suitable incremental step size in nonlinear analysis, particularly for plasticity problems, (c) convergence to benchmarks requires fine mesh since constant strain across first order VEM polygons makes fully plastic cross-sections difficult to achieve, (d) stabilization stiffness does seem to cause some interference to consistent linearization of the tangent stiffness matrix in plasticity problems. 

Suggested future work includes the following:  (a) investigate the effect of stabilization free VEM, (b) incorporate finite strains, (c) explore consistent linearization with projectors, (d) include problems of near incompressibility.

\section*{Acknowledgements}
\ack{
LLY acknowledges the research support of Walla Walla University.
Helpful discussions with N. Sukumar are also gratefully acknowledged.
}

\bibliography{bibCoVEM}

%
%
%
\end{document}